\documentclass[a4paper, 11pt]{amsart}
\usepackage{amsfonts}
\usepackage{amsfonts}
\usepackage{amsfonts}
\usepackage{amsfonts}
\usepackage{mathrsfs}
\usepackage{mathrsfs}
\usepackage{amsfonts,latexsym,rawfonts,amsmath,amssymb,amsthm,titletoc,appendix}
\usepackage[alphabetic]{amsrefs}
\usepackage{stmaryrd}  
\usepackage[applemac]{inputenc}
\usepackage[T1]{fontenc}
\usepackage{textcomp}
\usepackage{geometry}
\usepackage{graphicx}
\usepackage{mathptmx}
\usepackage{hyperref}
\usepackage[marginal]{footmisc}

\numberwithin{equation}{section}

\newcommand{\beq}{\begin{equation}}
\newcommand{\eeq}{\end{equation}}
\newcommand{\beqs}{\begin{eqnarray*}}
\newcommand{\eeqs}{\end{eqnarray*}}
\newcommand{\beqn}{\begin{eqnarray}}
\newcommand{\eeqn}{\end{eqnarray}}
\newcommand{\beqa}{\begin{array}}
\newcommand{\eeqa}{\end{array}}

\def\diag{{\rm diag}}

\newtheorem{prop}{Proposition}[section]
\newtheorem{theo}[prop]{Theorem}
\newtheorem{lem}[prop]{Lemma}

\begin{document}
\title{New area-minimizing Lawson-Osserman cones}
\author{Xiaowei Xu \quad Ling Yang\quad Yongsheng Zhang $\dagger$}\thanks{$\dagger$ Y. Zhang is the corresponding author.}
\address{
School of Mathematical Sciences,  University of Science and Technology of China, Hefei, 230026, Anhui
province, China; and Wu Wen-Tsun Key Laboratory of Mathematics, USTC, Chinese Academy of
Sciences, Hefei, 230026, Anhui, P.R. China.}\email
{
xwxu09@ustc.edu.cn.
}
\address{
School of Mathematical Sciences, Fudan University,
Shanghai 200433, P.R. China.}\email
{
yanglingfd@fudan.edu.cn}
\address{
 School of Mathematics and Statistics, Northeast Normal University, Changchun 130024, Jilin province, P.R. China.}\email
 {
yongsheng.chang@gmail.com}


\maketitle

\noindent \textbf{Abstract:}
It has been 40 years since 
Lawson and Osserman introduced the three minimal cones associated with Dirichlet problems
in their 1977 Acta paper \cite{[LO77]}.
The first cone was shown area-minimizing by Harvey and Lawson in the celebrated paper \cite{[HL82]}.
In this paper, we confirm that the other two are also area-minimizing.
In fact, we show that every Lawson-Osserman cone of type $(n,p,2)$ constructed in \cite{[XYZ16]} is area-minimizing.

\vspace{0.2cm}

\vspace{0.2cm}


\tableofcontents
\titlecontents{section}[0em]{}{\hspace{.5em}}{}{\titlerule*[1pc]{.}\contentspage}
\titlecontents{subsection}[1.5em]{}{\hspace{.5em}}{}{\titlerule*[1pc]{.}\contentspage}

\section{Introduction}

Let $\Sigma\subset S^N\subset \textbf{R}^{N+1}$ be an oriented closed submanifold (or rectifiable current) in the unit sphere.
Then 
 \[
                  \mathcal{C}\Sigma=\{tx:t\in[0,\infty) \text{ and } x\in \Sigma\}
                  \]
                  is called the cone over $\Sigma$.
                  We say $\mathcal{C}\Sigma$ is area-minimizing,
                  if the truncated cone inside the unit ball has least area
                  among all integral currents with boundary $\Sigma$.

       The study of area-minimizing cones is a central topic in the geometric measure theory.
        By the well-known result of Federer
        (Theorem 5.4.3 in \cite{F},
        also see Theorem 35.1 and Remark 34.6 (2) in Simon \cite{LS})
        that a tangent cone at a point of an area-minimizing rectifiable current is itself area-minimizing,
        it is meaningful to explore the diversity of area-minimizing cones
        for better understandings about local behaviors of 
        area-minimizing integral currents.

      Area-minimizing cones also capture
      behaviors at infinity for area-minimizing surfaces in certain cases.
      The celebrated Bernstein problem stimulates the study on
      the nonexistence, the existence and the diversity of area-minimizing hypercones,
      e.g.
       \cite{Fleming, DeG, Almgren, Simons, BDG, BL, PS1, PS2, HS, FK, LS89, [L91]}.
       In contrast, not
       {quite}
       much is known about area-minimizing cones of higher codimension{s}.
       Following \cite{BL}, Cheng \cite{Cheng} found homogeneous area-minimizing cones of codimension two.
       Around the same time, Lawlor \cite{[L91]} developed a systematic method,
     called the curvature criterion,
     for determining
     whether a minimal cone is indeed area-minimizing, 
       for instance, the classification of area-minimizing cones over products of spheres
       and the first examples of minimizing cones over nonorientable surfaces in the sense of mod 2.

       Among others,
       three interesting non-parametric minimal cones were
       constructed by Lawson and Osserman \cite{[LO77]} as follows.
       Let $\eta,\;\eta'$ and $\eta''$ denote the (normalized) Hopf maps
       $S^{2^i-1}$ $\rightarrow$ $S^{2^{i-1}}$ for $i=2,\;3 $ and $4$.
      Then Lawson and Osserman considered 
      the minimal embeddings
             \begin{equation}\label{LO3C}
             S^{2^i-1}\rightarrow S^{2^i+2^{i-1}},\ \ \ \ x \mapsto \Big(\alpha_2 x, \beta_2 \eta(x)\Big),\;\;
             \Big(\alpha_3 x, \beta_3 \eta'(x)\Big),\;\;\Big(\alpha_4 x, \beta_4 \eta''(x)\Big)
             \end{equation}
             with $\alpha_i=\sqrt{\frac{4(2^{i-1}-1)}{3(2^i-1)}}$ and $\beta_i=\sqrt{\frac{2^i+1}{3(2^i-1)}}$
       respectively.
       Over these minimal spheres, three minimal cones $C_i$ for $i=2,\ 3$ and $4$ are then obtained. 
       They, respectively, produce Lipschitz (not $C^1$) solutions to the Dirichlet problems of minimal surfaces
       over unit disks for boundary data
             \[
             \phi=\frac{\beta_2}{\alpha_2}\eta,\ \ \ \ \phi'=\frac{\beta_3}{\alpha_3}\eta',\ \ \ \phi''=\frac{\beta_4}{\alpha_4}\eta''.
             \]
       It was shown later in \cite{[HL82]} that $C_2$ is calibrated
       by the so-called coassociative calibration and, hence, area-minimizing in $\textbf{R}^7$ by the fundamental theorem of calibrated geometry.
       This canonical calibration in some way exhibits a special interaction between algebraic and geometric structures of $\textbf{R}^7$.
       Due to the lack of similar understandings,
       it remained open
       for decades
       whether the other two cones are area-minimizing.
       In this paper, we answer both affirmatively.
       
       Inspired by \cite{[LO77]},
     uncountably many non-parametric minimal cones of Lawson-Osserman type
       are constructed in \cite{[XYZ16]}.
       Each of them is a twisted graph similar to those in \eqref{LO3C}
       for the composition of a Hopf fibration and
      a homothetic
      (i.e., up to a constant factor, isometric)
      minimal
       immersion of
      even degree
       from projective space into a unit sphere.
       In the present paper we focus on,
       in particular,
       the Lawson-Osserman cones derived from standard minimal immersions of degree 2.
       More explicitly, let $\pi^{2n+1,2n}$ be the Hopf fibration from $S^{2n+1}$ onto $(\textbf{CP}^n,g_{FS})$
       and $\Phi_s$ the standard minimal immersion from $(\textbf{CP}^n,g_{FS})$ into a sphere.
        Then there is a unique positive rescaling $f$
         of $\Phi_s\circ \pi^{2n+1,2n}$,
         which maps
        $S^{2n+1}$ into $S^{n^2+2n-1}$.
        In fact, $f$ is exactly $\Phi\circ\pi^{2n+1,2n}$ in our \S3,
        each component of which is an $S^1$-invariant
        harmonic polynomial of degree 2,
        and the associated Lawson-Osserman sphere is
       $$F:S^{2n+1}\longrightarrow S^{(n+1)^2+2n},
           \hspace{0.5cm}
           z\mapsto\Big(a_n z,\;\;
           b_n f(z) \Big),$$
      where $a_n=\sqrt{\frac{2(n+1)}{(2n+1)(n+2)}}$ and $b_n=\sqrt{\frac{n(2n+3)}{(2n+1)(n+2)}}$
      (see Theorem 2.3 of \cite{[XYZ16]}).
      Similarly,
      we get 
      Lawson-Osserman spheres
      $F':S^{4n+3}\longrightarrow S^{2n^2+7n+3}$ and $F'':S^{15}\longrightarrow S^{24}$
      associated to Hopf fibrations $\pi^{4n+3,4n}: S^{4n+3}\longrightarrow \textbf{HP}^n$ and $\pi^{15,8}:S^{15}\longrightarrow \textbf{OP}^1$ respectively.
     {
     We remark that
      every isometric minimal immersion of degree 2
      from a projective space
      (endowed with certain multiple of the standard metric)
      into a unit sphere 
      is congruent to a standard one by an isometry of the target sphere.
      See \cite{[DW71], [M80], [M81], [U85], O} for details.
     Hence, we only need to consider the standard cases above.
     }
      Let us denote 
       by $\mathcal{C}F$, $\mathcal{C}F'$ and $\mathcal{C}F''$
      the Lawson-Osserman cones
     over
      the images of $F$, $F'$ and $F''$
     {respectively}.
     It is worth noticing that
          $\mathcal{C}F$, $\mathcal{C}F'$ for $n=1$ and $\mathcal{C}F''$ are exactly
          $C_2$, $C_3$ and $C_4$
     constructed by Lawson and Osserman. 
     In this paper, we establish  

     \vspace{0.3cm}
     \textbf{Main Theorem}. \emph{The minimal cones $\mathcal{C}F$, $\mathcal{C}F'$ and $\mathcal{C}F''$ are area-minimizing.}
       \vspace{0.3cm}

     The understanding of geometries
     of $F$, $F'$ and $F''$ is a key step in the proof of our main theorem.
     Since they are all homogeneous according to the classical theory of representation,
     we do calculations only at a base point.
     Via a good choice of local parameterizations of odd-dimensional spheres,
     we get a local adapted frame. 
     Then, by taking second derivatives,
     we gain the second fundamental form. 
     This method works for all the three types of Lawson-Osserman
     spheres uniformly, and seems more readable,
     for our purpose,
     than calculations
     through the theory of Lie group and Lie algebra.

       The paper is organized as follows.
       In \S 2 we introduce preferred local parameterizations of odd-dimensional spheres and
       make agreements on notations.
In \S 3, three types of Lawson-Osserman spheres are constructed explicitly and their second fundamental forms are calculated.
       By explanations on Lawlor's curvature criterion in \S 4.1 and the computation results in \S 3,
       we show in \S 4.2 that
       vanishing angles exist for 
       $\mathcal{C}F$, $\mathcal{C}F'$ and $\mathcal{C}F''$,
       with the exception of $C_2$, 
       and
       that the corresponding normal wedges do not intersect.
       Hence, by \cite{[HL82]} and Lawlor's curvature criterion,
       we accomplish our main theorem.
{\ }
\\

\section{Quaternions, octonions and odd-dimensional spheres}\label{sec2}
         We recall some basic facts about quaternions and octonions
         and make a choice of preferred local parameterizations
         of odd-dimensional spheres at a fixed point by exponential map.

            Let $\textbf{R}$, $\textbf{C}$ be the real and complex number fields, and
            $\textbf{R}^n$, $\textbf{C}^n$ the real and complex $n$-tuple spaces respectively.
            Conventionally, we identify $\textbf{C}^n$ with $\textbf{R}^{2n}$ by $(z_1,\ldots,z_n)$ $\mapsto$ $(x_1,\ldots,$ $x_n,x_{n+1},$ $\ldots,x_{2n})$,
            where $z_k=x_k+\texttt{i}x_{n+k}\in \textbf{C}$ and $\texttt{i}^2=-1$.
           Let $\textbf{H}$ be the real division algebra of quaternions.
           An element of $\textbf{H}$ can be written uniquely as $a=z_1+z_2\texttt{j}$,
           where $z_1,z_2\in\textbf{C}$ and $\texttt{j}\in\textbf{H}$ satisfies
           $\texttt{j}^2=-1,\;\;z\texttt{j}=\texttt{j}\bar{z},\ $
           for all $z\in\textbf{C}$.
           In this way,
                   $\textbf{C}$ sits inside $\textbf{H}$ as a subalgebra and
                   $\textbf{H}$ becomes a complex vector space
                   by $\textbf{C}$-action on the left.
                   Thus, the $n$-tuples $\textbf{H}^n$ can be identified with $\textbf{C}^{2n}$ by $(a_1,\ldots,a_n)\mapsto(z_1,z_2\ldots,z_{2n-1},z_{2n})$, where $a_k=z_{2k-1}+z_{2k}\texttt{j}\in\textbf{H}$. Let $\textbf{O}$ be the algebra of octonions.
                    Then elements in $\textbf{O}$ can be written as $p=a_1+a_2\texttt{e}$, where $a_1,a_2\in\textbf{H}$, $\texttt{e}\in\textbf{O}$ and $\texttt{e}^2=-1$, $p\texttt{e}=-a_2+a_1\texttt{e}$. For another $q=a_3+a_4\texttt{e}\in\textbf{O}$, we have $pq=(a_1a_3-\bar{a}_4a_2)+(a_4a_1+a_2\bar{a}_3)\texttt{e}$. Similarly, the $n$-tuples $\textbf{O}^n$ can be identified with $\textbf{H}^{2n}$ by $(p_1,\ldots,p_n)\mapsto(a_1,a_2\ldots,a_{2n-1},a_{2n})$, where $p_k=a_{2k-1}+a_{2k}\texttt{e}\in\textbf{O}$.
The conjugate of a quaternion (resp. octonion) is defined 
                to be $\bar a=\bar{z}_1-z_2\texttt{j}$ (resp. $\bar p=\bar{a}_1-a_2\texttt{e}$).
                Hence, the norm square of a quaternion (resp. octonion) is given by $|a|^2=a\bar a$ (resp. $|p|^2=p\bar p$).

           An odd-dimensional unit sphere can be described by
\begin{equation*}
            S^{2m+1}(1)=\Big\{(z_1,\ldots,z_{m+1})\in \textbf{C}^{m+1}\;\Big|\;\;\;\sum\limits_{k=1}^{m+1}|z_k|^2=1\Big\}.
\end{equation*}
             For a local parameterization at point $p=(1,0,\ldots,0)\in S^{2m+1}$, we set
\begin{equation*}
X_1=\texttt{i}(\mbox{E}_{1\;1}-\mbox{E}_{2\;2}),\;\;\;\;X_k=\mbox{E}_{1\;k}-\mbox{E}_{k\;1},\;\;\;\;X_{m+k}=\texttt{i}(\mbox{E}_{1\;k}+\mbox{E}_{k\;1}),
\end{equation*}
where $2\leq k\leq m+1$, $\mbox{E}_{k\;l}$ is the $(m+1)\times(m+1)$-matrix
                 with value 1 in the $(k,l)$-slot and zero for others.
Then, it is easy to check that
\begin{equation*}
              p(t)=p\;e^{t_1X_1}e^{t_2X_2}\cdots e^{t_{2m+1}X_{2m+1}}\ \
              \text{ for small } t=(t_1,\ldots,t_{2m+1})\in \textbf{R}^{2m+1}
\end{equation*}
          parameterizes {a} neighborhood of $p$.
          Writing $p(t)=(z_1,\ldots,z_{m+1})$, through direct calculations,
          we have
\begin{eqnarray}
z_1&=&\big(1-\frac{1}{2}\sum\limits_{A=1}^{2m+1}t^2_A\big)+\texttt{i}\big(t_1+\sum\limits_{k=2}^{m+1}t_k\,t_{m+k}\big)+o(t^2),\\
z_k&=&\big(t_k-t_1t_{m+k}\big)+\texttt{i}\big(t_{m+k}+t_1t_k\big)+o(t^2),
\end{eqnarray}
where $2\leq k\leq m+1$ and $t^2=\sum\limits_{A=1}^{2m+1}t^2_A$.
          In our parameterization, the tangent space $T_pS^{2m+1}$ is spanned by
           $\{\epsilon_A=\frac{\partial p(t)}{\partial t_A}\Big|_{t=0}\}$, more explicitly, by
\begin{equation*}
\epsilon_1=\texttt{i}\mbox{E}_1\in\textbf{C}^{m+1},\hspace{1cm}\mbox{   or   }\hspace{1cm} (0,\mbox{E}_1)\in\textbf{R}^{2m+2},
\end{equation*}
\begin{equation*}
\epsilon_k=\mbox{E}_k\in\textbf{C}^{m+1},\hspace{1.2cm}\mbox{   or   }\hspace{1cm} (\mbox{E}_k,0)\in\textbf{R}^{2m+2},
\end{equation*}
 \begin{equation*}
\epsilon_{m+k}=\texttt{i}\mbox{E}_k\in\textbf{C}^{m+1},\hspace{0.6cm}\mbox{   or   }\hspace{1cm} (0,\mbox{E}_k)\in\textbf{R}^{2m+2},
\end{equation*}
where $2\leq k\leq m+1$ and $\mbox{E}_k$ is the vector in $\textbf{R}^{m+1}$
             with value 1 in the $k$-th position and zero for others.
{\ }\\{\ }\\

\section{On explicit construction of Lawson-Osserman spheres}\label{sec3}

              In this section we shall give more details in what we are concerned with.
              As explained in the introduction,
              embedded minimal spheres can be built
              on the compositions of Hopf fibrations and
              minimal immersions from
              complex projective spaces,
              quaternion projective spaces
              or the Cayley projective line into unit spheres.
              Since
              the thought originates from Lawson and Osserman's constructions,
              we call such minimal spheres
               \textbf{Lawson-Osserman spheres} and 
              the associated cones 
              \textbf{Lawson-Osserman cones}.
              By the fact that there are only three families of Hopf fibrations
              $\pi^{2n+1,2n}: S^{2n+1}\longrightarrow \textbf{CP}^n$,
              $\pi^{4n+3,4n}: S^{4n+3}\longrightarrow \textbf{HP}^n$ and
              $\pi^{15,8}:S^{15}\longrightarrow \textbf{OP}^1$,
              we, accordingly, divide
               Lawson-Osserman spheres into three types.
              It is worth mentioning
              that Tang \cite{[T01]} proved
              {the nonexistence of submersion}
              from $S^{23}$ to the Cayley projective plane.

             Let us first recall the classical Lawson-Osserman construction.
             By the Hopf map $\eta:S^3\longrightarrow S^2$,
             $(z_1,z_2)\mapsto(|z_1|^2-|z_2|^2,2\bar z_1 z_2)$,
             Lawson and Osserman \cite{[LO77]} gave a minimal immersion
             \begin{equation*}
             F:S^3(1)\longrightarrow S^6(1),\;\;\;z\mapsto\big(\frac{2}{3}z,\frac{\sqrt{5}}{3}\eta(z)\big).
             \end{equation*}
             Later, in \cite{[HL82]}, Harvey and Lawson proved that the cone of $F$ is area-minimizing.
             Therefore, the topological space $\textbf{R}^4$ can emerge in $\textbf{R}^7$
             as a nontrivial area-minimizing cone.
             To find more such minimal immersions,
             we observe $\eta$ from another perspective.
             In fact, it is a composition of the Hopf fibration $\pi^{3,2}$ and
             a degree 2 map $\Phi$ from $\textbf{CP}^1$ into $S^2$.
             Explicitly,
             \begin{equation*}
             \pi^{3,2}:S^3\longrightarrow \textbf{CP}^1,\;(z_1,z_2)\mapsto[z_1,z_2],
             \;\;\text{and}\;\;\Phi:\textbf{CP}^1\longrightarrow S^2,\;[z_1,z_2]\mapsto(|z_1|^2-|z_2|^2,2\bar z_1 z_2).
             \end{equation*}
             Analogous constructions in \cite{[LO77]}
             were also given for Hopf maps
             $\eta':S^7\longrightarrow S^4$,
             $(a_1,a_2)\mapsto(|a_1|^2-|a_2|^2,2\bar a_1 a_2)$ and
             $\eta'':S^{15}\longrightarrow S^8$,
             $(p_1,p_2)\mapsto(|p_1|^2-|p_2|^2,2\bar p_1 p_2)$.
             However, to our knowledge,
             it was unknown whether
             the cones associated to $\eta'$ and $\eta''$ are area-minimizing.
             By results to be established in this section,
             we can have positive conclusions for both in \S \ref{sec4}.

             By composing Hopf fibrations and
             {homothetic} minimal immersions of degree 2 (explained below)
             from $\textbf{CP}^n$, $\textbf{HP}^n$ 
             into unit spheres,
             we gain lots of Lawson-Osserman spheres.
             The minimality follows from a general theorem of authors \cite{[XYZ16]}.
             Notice that $\textbf{CP}^n$, $\textbf{HP}^n$ and unit spheres
             are compact symmetric spaces.
             Hence, such immersions can be realized by equivariant ones
             compatible with their Lie group structures.

             Now, we review standard minimal immersions of
             a compact irreducible Riemannian symmetric space $(M,g)$
             into unit spheres (see \cite{[DW71], W, [U85]}).
             Let $\Delta$ be the Laplace-Beltrami operator of $(M,g)$ acting on $C^\infty$-functions,
             $\lambda_k$ the $k$-th eigenvalue of $\Delta$ with $0=\lambda_0<\lambda_1<\cdots$,
             and $V^k$ the corresponding eigenspace.
             Set $\mbox{dim}V^k=m(k)+1$.
             Then one can define an $L^2$-inner product on
             $V^k$ by
                    \[(f,h):=\int_{M}\;fh\;d\mu\]
              where $d\mu$ of $(M,g)$ is the normalized canonical measure
              with $\int_{M}d\mu=m(k)+1$.
              Suppose $\{f_0,\ldots,f_{m(k)}\}$ form an orthonormal basis of $V^k$.
              By Takahashi's Theorem \cite{[T66]},
              the standard map $x_k$ from $M$ to $\textbf{R}^{m(k)+1}$
obtained              by sending $p\in M$ to $(f_0(p),\ldots,f_{m(k)}(p))$
              gives an isometric minimal immersion of $(M,\frac{\lambda_k}{dim\, M}\;g)$ into $S^{m(k)}(1)$.
              This standard minimal immersion $x_k$ can also be understood as follows.
              Let $(G,K)$ be a symmetric pair with $M=G/K$.
              Then, a point of $M$ can be regarded as $\sigma K$ for some $\sigma\in G$,
              and $G$ acts on $V^k$ by $(\sigma\cdot f)(p)=f(\sigma^{-1}p)$, $\sigma\in G$, $p\in M$.
              In this way,
              an orthogonal representation of $G$
              is given in terms of $\sigma\cdot f_\alpha=\sum\limits_{\alpha=0}^{m(k)}a_{\alpha\beta}f_\beta$
               by
               \begin{equation}\label{3.1}
               \rho_k:G\longrightarrow O(m(k)+1),\hspace{0.5cm}\sigma\mapsto (a_{\alpha\,\beta}(\sigma)).
               \end{equation}
               Up to a rigidity of $V^k$,
               one can assume that $\big(f_1(eK),\ldots,f_{m(k)}(eK)\big)=\mbox{E}_1=\big(1,0,\ldots,0\big)$.
               Then we get $x_k(\sigma K)=\mbox{E}_1\rho_k(\sigma)$, $\sigma\in G$.

               We give an alternative description on $V^k$
               for $M=\textbf{CP}^n$ equipped with the Fubini-Study metric.
               Let $\phi$ be a complex valued homogeneous polynomial
               in $2n+2$ variables $z_1,\ldots,z_{n+1},\bar{z}_1,\ldots,\\\bar{z}_{n+1}$
               on $\textbf{C}^{n+1}$.
               It is said to be of $(p,q)$-type if
               $$\phi(cz_1,\ldots,cz_{n+1},\bar{c}\bar{z}_1,\ldots,\bar{c}\bar{z}_{n+1})=
               c^p\bar{c}^q\phi(z_1,\ldots,z_{n+1},\bar{z}_1,\ldots,\bar{z}_{n+1}),
               \text{ for } \forall\; c\in\textbf{C}.$$
               Denote by $P^{n+1}_{p,q}$ the complex vector space of
               all homogeneous polynomials of $(p,q)$-type on $\textbf{C}^{n+1}$.
               Note that functions in $P^{n+1}_{q,q}$ are $S^1$-invariant.
               So, they descend to functions on $\textbf{CP}^n$.
               Set
               $$
               H^{n+1}_{p,q}=\big\{\phi\in P^{n+1}_{p,q}\;\big|\;D\phi=0\big\},\ \text{ where }
               D=-4\sum\limits_{k=1}^{n+1}\frac{\partial^2}{\partial z_k\bar z_k}.
               $$
                It is a well-known fact (see \cite{[M80], [U85]})
                that the $k$-th eigenspace $V^k$ is $SU(n+1)$-isomorphism to
                $H^{n+1}_{k,k}\cap C^\infty(\textbf{C}^{n+1},\textbf{R})$ for $\textbf{CP}^n$
                via the Hopf fibration $\pi^{2n+1,2n}$,
                where $C^\infty(\textbf{C}^{n+1},\textbf{R})$
                consists of all real valued $C^\infty$-functions on $\textbf{C}^{n+1}$.
                 In the present paper, we shall focus on $H^{n+1}_{1,1}\cap C^\infty(\textbf{C}^{n+1},\textbf{R})$,
                the space of isometric minimal immersions of degree 2.
                There is also a similar description in terms of quaternions
                valid for $\textbf{HP}^n$ equipped with the canonical metric.

               Throughout this paper, we will, if not otherwise specified,
               use the following convention for indices:
\begin{equation*}
1\leq A,B,\cdots\leq 2m+1;\;\;\;2\leq k,l,\cdots\leq n+1;\;\;\;2\leq\alpha,\beta,\cdots\leq n,
\end{equation*}
where $m$ will take $n$ or $2n+1$ in the sequel.
{\ }\\

\subsection{Type-I Lawson-Osserman spheres}\label{s3.1}
           To obtain a standard minimal immersion from $\textbf{CP}^n$ into unit sphere,
           we need to find an orthonormal basis of $H^{n+1}_{1,1}\cap C^\infty(\textbf{C}^{n+1},\textbf{R})$.
           For an element $\phi\in H^{n+1}_{1,1}\cap C^\infty(\textbf{C}^{n+1},\textbf{R})$,
           it follows that
           $\phi=\sum\limits_{k,l=1}^{n+1}\lambda_{k\bar l}\;z_k\bar z_l$
           subject to
           $D\phi=0$ and $\bar\phi=\phi$.
           Set $\lambda_{\bar l k}=\bar \lambda_{l\bar k}$.
           The requirements become
           \begin{equation}\label{3.2}
           \sum\limits_{k=1}^{n+1}\lambda_{k\bar k}=0,\hspace{0.5cm} \lambda_{k\bar l}=\lambda_{\bar l k}.
           \end{equation}
           With this
           {comprehension}, we have:
           \begin{lem}\label{lemma3.1}
           There is an orthonormal basis of $H^{n+1}_{1,1}\cap C^\infty(\textbf{C}^{n+1},\textbf{R})$ w.r.t. the $L^2$-inner product given by
\begin{equation*}
\phi_\alpha=c_{n,\alpha}\Big(|z_\alpha|^2-\frac{1}{n\!+\!1\!-\!\alpha}\;\sum\limits_{k=\alpha+1}^{n+1}|z_k|^2\Big),\hspace{0.5cm}1\leq \alpha\leq n,
\end{equation*}
\begin{equation*}
\phi_{k\,l}=d_n\;\mbox{\emph{Re}}(z_k\,\bar z_l),\hspace{0.5cm} \phi_{\bar{k}\,\bar{l}}=d_n\;\mbox{\emph{Im}}(z_k\,\bar z_l),\hspace{0.5cm}1\leq k<l\leq n+1,
\end{equation*}
where $c_{n,\alpha}=\sqrt{\frac{(n+1)(n+1-\alpha)}{n(n+2-\alpha)}}$ and $d_n=\sqrt{\frac{2(n+1)}{n}}$.
\end{lem}

\emph{Proof}. From (\ref{3.2}), we know that $u_\alpha=|z_\alpha|^2-|z_{\alpha+1}|^2$, $v_{k\,l}=\mbox{Re}(z_k\,\bar z_l)$ and $v_{\bar{k}\,\bar{l}}=\mbox{Im}(z_k\,\bar z_l)$ form a basis of $H^{n+1}_{1,1}\cap C^\infty(\textbf{C}^{n+1},\textbf{R})$. Then, following the Schmidt orthonormalization process w.r.t. the $L^2$-inner product, one can get
           $\{\phi_\alpha,\;\phi_{k\,l},\;\phi_{\bar k\,\bar l}\}$
           as an orthonormal basis of $H^{n+1}_{1,1}\cap C^\infty(\textbf{C}^{n+1},\textbf{R})$.
           We leave details to readers.\hfill$\Box$

\vspace{0.3cm}

           Hence, we gain an isometric minimal immersion
           $\Phi:\textbf{CP}^n\longrightarrow S^{n(n+2)-1}(1)$ expressed by
           $[z]\mapsto\Big(\phi_\alpha(z),\; \phi_{k\,l}(z),\;\phi_{\bar k\,\bar l}(z)\Big)$
           and a Lawson-Osserman sphere given by
            $f=\Phi\circ\pi^{2n+1,2n}:S^{2n+1}(1)$ $\rightarrow$ $S^{n(n+2)-1}(1)$.
           In particular,
           $f$ is just the Hopf map $\eta$ when $n=1$.
           The type-I Lawson-Osserman sphere is represented by
           \begin{equation*}
           F:S^{2n+1}(1)\longrightarrow S^{(n+1)^2+2n}(1),
           \hspace{0.5cm}
           z\mapsto\Big(a_n\big(\mbox{Re}(z),\mbox{Im}(z)\big),\;\;
           b_n\big(\phi_\alpha(z),\; \phi_{k\,l}(z),\;\phi_{\bar k\,\bar l}(z)\big)\Big),
           \end{equation*}
           where $a_n=\sqrt{\frac{2(n+1)}{(2n+1)(n+2)}}$ and $b_n=\sqrt{\frac{n(2n+3)}{(2n+1)(n+2)}}$.

\emph{Remarks}. (a) For coefficients $a_n$ and $b_n$
            which uniquely determine the minimality
            we refer to our recent paper \cite{[XYZ16]}.
            We will just check the minimality property in Proposition \ref{3.2};

             (b) By the construction of $\Phi$,
              $F$ is in fact homogeneous and
               its image is an $\mbox{Id}\oplus\rho_1\big(SU(n+1)\big)$-orbit through the base point
               $$P=\Big(a_n\,\mbox{E}_1,\;0,\;b_n\mbox{E}_1,\;0,\;0\Big),$$
               where $\rho_1$ is defined in (\ref{3.1}).
               Therefore, we will study its geometry only around one point for our purpose.

\vspace{0.3cm}

            More precisely, we shall compute the second fundamental form of $F$ at $P$
            in the remaining part of this subsection.
            Substituting (2.1) and (2.2) for $m=n$ into expressions of
            $\phi_\alpha,\;\phi_{k\,l},\;\phi_{\bar k\,\bar l}$, we obtain
\begin{eqnarray}\label{3.3}
\phi_1&=&1-\frac{n+1}{n}\,\sum\limits_{A=2}^{2n+1}t^2_A+o(t^2),\nonumber\\
\phi_\alpha&=&c_{n,\alpha}\Big[t_\alpha^2+t^2_{n+\alpha}-\frac{1}{n\!+\!1\!-\!\alpha}\sum\limits_{k=\alpha+1}^{n+1}\big(t_k^2+t_{n+k}^2\big)\Big]+o(t^2),
\nonumber\\
\phi_{1\,k}&=&d_n\,t_k+o(t^2),\hspace{0.3cm}\phi_{1\,\bar k}=-d_n\,t_{n+k}+o(t^2),\\
\phi_{k\,l}&=&d_n\Big(t_k\,t_l+t_{n+k}\,t_{n+l}\Big)+o(t^2),\nonumber\\
\phi_{\bar{k}\,\bar{l}}&=&d_n\Big(-t_k\,t_{n+l}+t_{l}\,t_{n+k}\Big)+o(t^2),\nonumber
\end{eqnarray}
where $2\leq\alpha\leq n$, $2\leq k\leq n+1$ and $2\leq k<l\leq n+1$.
             Noticing that $F_*(\epsilon_A)=\frac{\partial F}{\partial t_A}\Big|_{t=0}$,
             by (\ref{3.3}), (2.1) and (2.2), we have
\begin{eqnarray*}
F_*(\epsilon_1)&=&\Big(0,\;a_n\,\mbox{E}_1,\;0,\;0,\;0\Big),\hspace{0.3cm}\\
F_*(\epsilon_k)&=&\Big(a_n\,\mbox{E}_k,\;0,\;0,\;b_nd_n\,\mbox{E}_{1\;k},\;0\Big),\hspace{0.3cm}\\
F_*(\epsilon_{n+k})&=&\Big(0,\;a_n\,\mbox{E}_k,\;0,\;0,\;-b_nd_n\,\mbox{E}_{1\;k}\Big),
\end{eqnarray*}
where $2\leq k\leq n+1$.
             By normalization, we get an orthonormal basis of $F_*(T_pS^{2n+1})$:
\begin{equation*}
e_1=\frac{1}{a_n}\,F_*(\epsilon_1),\hspace{0.3cm}e_k=\frac{1}{\sqrt{a_n^2+b_n^2d_n^2}}\,F_*(\epsilon_k),
\hspace{0.3cm}e_{n+k}=\frac{1}{\sqrt{a_n^2+b_n^2d_n^2}}\,F_*(\epsilon_{n+k});
\end{equation*}
               and an orthonormal basis for the normal space of $F_*(T_pS^{2n+1})$ in $T_PS^{(n+1)^2+2n}$:
\begin{eqnarray}\label{3.3-1}
e_{2n+2}&=&\Big(-b_n\,\mbox{E}_1,\;0,\;a_n\,\mbox{E}_1,\;0,\;0\Big),\nonumber\\
e_{2n+1+k}&=&\frac{1}{\sqrt{a_n^2+b_n^2d_n^2}}\Big(-b_nd_n\,\mbox{E}_k,\;0,\;0,
\;a_n\,\mbox{E}_{1\,k},\;0\Big),\nonumber\\
e_{3n+1+k}&=&\frac{1}{\sqrt{a_n^2+b_n^2d_n^2}}\,\Big(0,\;b_nd_n\,\mbox{E}_{k},\;0,\;0,\;a_n\,\mbox{E}_{1\,k}\Big),\\
e_{4n+1+\alpha}&=&\Big(0,\;0,\;\mbox{E}_\alpha,\;0,\;0\Big),\nonumber\\
e_{k\,l}&=&\Big(0,\;\;0,\;\;0,\;\;\mbox{E}_{k\,l},\;\;0\Big),\nonumber\\
e_{\bar{k}\,\bar{l}}&=&\Big(0,\;\;0,\;\;0,\;\;0,\;\;\mbox{E}_{k\,l}\Big),\nonumber
\end{eqnarray}
where $2\leq k\leq n+1$, $2\leq\alpha\leq n$ and $2\leq k<l\leq n+1$.

\vspace{0.3cm}

             Set $F_{AB}=\frac{\partial^2F}{\partial t_A\partial t_B}\Big|_{t=0}$.
             Through direct computations, we obtain
\begin{eqnarray*}
F_{1\,1}&=&\Big(-a_n\,\mbox{E}_1,\;\;0,\;\;0,\;\;0,\;\;0\Big),\\
F_{\alpha\,\alpha}=F_{n+\alpha\;n+\alpha}&=&\Big(-a_n\,\mbox{E}_1,\;0,\;-\frac{2(n\!+\!1)b_n}{n}\mbox{E}_1-2b_n\sum\limits_{\beta=2}^{\alpha-1}
\frac{c_{n,\beta}}{n\!+\!1\!-\!\beta}\mbox{E}_\beta+2b_nc_{n,\alpha}\mbox{E}_\alpha,\;0,\;0\Big),\\
F_{n+1\,n+1}=F_{2n+1\;2n+1}&=&\Big(-a_n\,\mbox{E}_1,\;0,\;-\frac{2(n\!+\!1)b_n}{n}\mbox{E}_1-2b_n\sum\limits_{\beta=2}^{n-1}
\frac{c_{n,\beta}}{n\!+\!1\!-\!\beta}\mbox{E}_\beta-2b_nc_{n,n}\mbox{E}_n,\;0,\;0\Big),
\end{eqnarray*}
\begin{eqnarray*}
F_{1\,k}=\Big(0,\;\;a_n\mbox{E}_k,\;\;0,\;\;0,\;\;0\Big),\ \ \hspace{0.2cm}F_{1\,n+k}=\Big(-a_n\mbox{E}_k,\;\;0,\;\;0,\;\;0,\;\;0\Big),\ \
\hspace{0.2cm}F_{k\,n+k}=\Big(0,\;\;a_n\mbox{E}_1,\;\;0,\;\;0,\;\;0\Big),
\end{eqnarray*}
\begin{eqnarray*}
F_{k\,l}=F_{n+k\;n+l}=\Big(0,\;\;0,\;\;0,\;\;b_nd_n\mbox{E}_{k\,l},\;\;0\Big),\hspace{1.5cm}
F_{k\,n+l}=-F_{l\,n+k}=\Big(0,\;\;0,\;\;0,\;\;0,\;\;-b_nd_n\mbox{E}_{k\,l}\Big),
\end{eqnarray*}
where $2\leq\alpha\leq n$, $2\leq k\leq n+1$ and $2\leq k<l\leq n+1$. Define
\begin{equation*}
H_{AB}=\left\{
\begin{array}{lllll}
\frac{1}{a_n^2}F_{1\,1},\hspace{0.2cm}A=1,\;B=1\\
{}\\
\frac{1}{a_n\;\sqrt{a_n^2+b_n^2d_n^2}}F_{1\,B},\hspace{0.2cm}A=1,\;B\neq1,\\
{}\\
\frac{1}{a_n^2+b_n^2d_n^2}F_{A\,B},\hspace{0.2cm}A\neq1,\;B\neq1.
\end{array}
\right.
\end{equation*}
              Then, at $P$, the second fundamental form of $F$
              is given in terms of the frame
              $\{e_A,\;e_\tau,\;e_{k\,l},\;e_{\bar k\,\bar l}\;\;|\;
              1\leq A\leq 2n+1,\;2n+2\leq\tau\leq 5n+1,\;2\leq k<l\leq n+1\}$
               by
\begin{equation}\label{3.4}
h^\tau_{AB}=\big\langle H_{AB},e_\tau\big\rangle,\hspace{0.7cm}h^{(k,\,l)}_{AB}=\big\langle H_{AB},e_{k\,l}\big\rangle,\hspace{0.7cm} h^{(\bar k,\,\bar l)}_{AB}=\big\langle H_{AB},e_{\bar{k}\,\bar{l}}\big\rangle.
\end{equation}

\vspace{0.3cm}
More explicitly, we gain:

\begin{prop}\label{hforF}
The second fundamental form of $F$ at the base point $P$ w.r.t. the frame 
$\{e_A,\; e_\tau,\; e_{k\,l},\;\\ e_{\bar k\,\bar l}\;\;|\;1\leq A\leq 2n+1,\;2n+2\leq\tau\leq 5n+1,\;2\leq k<l\leq n+1\}$ is given by

\emph{(1)} \hspace{0.3cm} $h^{2n\!+\!2}_{1\,1}=\frac{b_n}{a_n}$,\hspace{0.3cm}   $h^{2n\!+\!2}_{AA}=-\frac{(n+2)a_nb_n}{n(a_n^2+b_n^2d_n^2)}$,\hspace{0.3cm}  $2\leq A\leq 2n+1;$

\emph{(2)} \hspace{0.3cm} $h^{2n\!+\!1\!+\!k}_{1\,n+k}=h^{3n\!+\!1\!+\!k}_{1\,k}=
\frac{b_nd_n}{a_n^2+b_n^2d_n^2}$,\hspace{0.3cm} $2\leq k\le n+1;$\hspace{0.3cm}

\emph{(3)} \hspace{0.3cm}$h^{4n\!+\!1\!+\!\alpha}_{k\,k}=h^{4n\!+\!1\!+\!\alpha}_{n+k\,n+k}=-\frac{2b_nc_{n,\alpha}}{(n\!+\!1-\alpha)(a_n^2+b_n^2d_n^2)}$,
\hspace{0.3cm} $2\leq \alpha\leq n$, \hspace{0.3cm}$\alpha<k\leq n+1;$

\emph{(4)} \hspace{0.3cm}$h^{4n\!+\!1\!+\!\alpha}_{\alpha\,\alpha}=h^{4n\!+\!1\!+\!\alpha}_{n+\alpha\,n+\alpha}=\frac{2b_nc_{n,\alpha}}{a_n^2+b_n^2d_n^2}$,
\hspace{0.3cm}$2\leq \alpha\leq n;$


\emph{(5)} \hspace{0.3cm}$h^{(k,l)}_{k\,l}=h^{(k,l)}_{n+k\,n+l}=h^{(\bar{k},\bar{l})}_{k\,n+l}=-h^{(\bar{k},\bar{l})}_{l\,n+k}
=\frac{b_nd_n}{a_n^2+b_n^2d_n^2},$\hspace{0.3cm} $2\leq k<l\leq n+1;$

\noindent  with the same value in the symmetric slot and zero for others.
\end{prop}

\emph{Proof}. It follows by direct computation.
                       Moreover, one can easily see that $F$ is minimal. \hfill$\Box$
{\ }\\

\subsection{Type-II Lawson-Osserman spheres}\label{s3.2}
                 In terms of quaternions,
                 \begin{equation*}
                 S^{4n+3}(1)=\Big\{(a_1,\ldots,a_{n+1})\in\textbf{H}^{n+1}\;\Big|\;\sum\limits_{k=1}^{n+1}|a_k|^2=1\Big\}.
\end{equation*}
                 A function in $H^{n+1}_{1,1}\cap C^\infty(\textbf{H}^{n+1},\textbf{R})$
                 restricted to $S^{4n+3}$ is $S^3$-invariant.
                 As a consequence, it descends to a function on $\textbf{HP}^n$.
                 Let $\phi\in H^{n+1}_{1,1}\cap C^\infty(\textbf{H}^{n+1},\textbf{R})$.
                 Then
                 $\phi=\sum\limits_{k,l=1}^{n+1}\lambda_{\bar k\, l}\;\bar a_k a_l$.
                 Notice that
                 $D=-4\sum\limits_{k=1}^{n+1}\frac{\partial^2}{\partial a_k\partial \bar a_k}$
                  in terms of quaternions,
                  so conditions $D\phi=0$ and $\bar\phi=\phi$ are equivalent to
\begin{equation*}
\sum\limits_{k=1}^{n+1}\lambda_{\bar k\, k}=0,\hspace{0.5cm} \lambda_{\bar k\, l}=\lambda_{l\,\bar k},
\end{equation*}
where $\lambda_{k\,\bar l}=\bar \lambda_{\bar k\, l}$. Similarly, we have:
\begin{lem}\label{lemma3.2}
There is an orthonormal basis of $H^{n+1}_{1,1}\cap C^\infty(\textbf{H}^{n+1},\textbf{R})$ w.r.t. the $L^2$-inner product given by
\begin{equation*}
\phi_\alpha=c_{n,\alpha}\Big(|z_{2\alpha-1}|^2+|z_{2\alpha}|^2-\frac{1}{n\!+\!1\!-\!\alpha}\;
\sum\limits_{k=\alpha+1}^{n+1}(|z_{2k-1}|^2+|z_{2k}|^2
{)}
\Big),\hspace{0.5cm}1\leq \alpha\leq n,
\end{equation*}
\begin{equation*}
\phi_{k\,l}=d_n\,\mbox{\emph{Re}}(\bar{z}_{2k-1}\bar z_{2l-1}+z_{2k}\bar z_{2l}),\hspace{0.3cm}
\phi_{\bar{k}\,\bar{l}}=d_n\,\mbox{\emph{Im}}(\bar{z}_{2k-1}\bar z_{2l-1}+z_{2k}\bar z_{2l}),\hspace{0.3cm}1\leq k<l\leq n+1,
\end{equation*}
\begin{equation*}
\tilde{\phi}_{k\,l}=d_n\,\mbox{\emph{Re}}(\bar{z}_{2k-1}\bar z_{2l}-z_{2k}\bar z_{2l-1}),\hspace{0.3cm}
\tilde{\phi}_{\bar{k}\,\bar{l}}=d_n\,\mbox{\emph{Im}}(\bar{z}_{2k-1}\bar z_{2l}-z_{2k}\bar z_{2l-1}),\hspace{0.3cm}1\leq k<l\leq n+1,
\end{equation*}
where $c_{n,\alpha}=\sqrt{\frac{(n+1)(n+1-\alpha)}{n(n+2-\alpha)}}$, $d_n=\sqrt{\frac{2(n+1)}{n}}$ and
                   we consider $a_k$ as $z_{2k-1}+z_{2k}\emph{\texttt{j}}$.
\end{lem}

\emph{Proof}. Similar to the proof of Lemma \ref{lemma3.1}. \hfill$\Box$

\vspace{0.3cm}

              Thus, we have an isometric minimal immersion
              $\Phi:\textbf{HP}^n\longrightarrow S^{2n^2+7n+3}(1)$ by
              $$[a]\mapsto\Big(\phi_\alpha(a),\; \phi_{k\,l}(a),\; \phi_{\bar k\,\bar l}(a),\;
              \tilde{\phi}_{k\,l}(a),\; \tilde{ \phi}_{\bar{k}\,\bar{l}}(a)\Big),$$
               and
            a Lawson-Osserman cone determined by $f'=\Phi\circ\pi^{4n+3,4n}$.
            It is known that $f'$ is just the Hopf map $\eta'$ when $n=1$.
            The type-II Lawson-Osserman sphere is
            represented by $F':S^{4n+3}\longrightarrow S^{2n^2+7n+3}$ sending
\begin{equation*}
a\mapsto\Big(\tilde{a}_n\big(\mbox{Re}
({a}),
\;\mbox{Im}
({a})
\big),\;\;\tilde{b}_n\big(\phi_\alpha(a),\; \phi_{k\,l}(a),\;\phi_{\bar k\,\bar l}(a),\;\tilde{\phi}_{k\,l}(a),\;\tilde{\phi}_{\bar k\,\bar l}(a)\big)\Big),
\end{equation*}
                where $\tilde{a}_n=\sqrt{\frac{6(n+1)}{(n+2)(4n+3)}}$ and $\tilde{b}_n=\sqrt{\frac{n(4n+5)}{(n+2)(4n+3)}}$.
                Again, for the choice of $\tilde{a}_n$ and $\tilde{b}_n$,
                we refer to \cite{[XYZ16]} for a general explanation.
                By the construction of $\Phi$,
                 $F$ is homogeneous and its image is an $\mbox{Id}\oplus\rho_1\big(Sp(n+1)\big)$-orbit
                 through the base point
                 $$P=\Big(a_n\,\mbox{E}_1,\;0,\;b_n\mbox{E}_1,\;0,\;0,\;0,\;0\Big),$$
                 where $\rho_1$ is defined in (\ref{3.1}).
                 Next, we compute its second fundamental form at $P$.

\vspace{0.3cm}

Using (2.1) and (2.2) with $m=2n+1$, we have
\begin{eqnarray}\label{3.6}
\phi_1&=&1-\frac{n+1}{n}\,\sum\limits_{k=2}^{n+1}\Big(t^2_{2k-1}+t^2_{2k}+t^2_{2n+2k}+t^2_{2n+1+2k}\Big)+o(t^2),\nonumber\\
\phi_\alpha&=&c_{n,\alpha}\Big[t_{2\alpha-1}^2+t^2_{2\alpha}+t^2_{2n+2\alpha}+t^2_{2n+1+2\alpha}-
\frac{1}{n\!+\!1\!-\!\alpha}\sum\limits_{k=\alpha+1}^{n+1}\big(t_{2k-1}^2+t^2_{2k}\nonumber\\
&{}&\hspace{1cm}+t^2_{2n+2k}+t^2_{2n+1+2k}\big)\Big]+o(t^2),
\nonumber\\
\phi_{1\,k}&=&d_n\Big(t_{2k-1}+t_2t_{2k}+t_{2n+3}t_{2n+1+2k}\Big)+o(t^2),\nonumber\\
\phi_{1\,\bar k}&=&d_n\Big(t_{2n+2k}-t_2t_{2n+1+2k}+t_{2k}t_{2n+3}\Big)+o(t^2),\nonumber\\
\phi_{k\,l}&=&d_n\Big(t_{2k-1}t_{2l-1}+t_{2k}t_{2l}+t_{2n+2k}t_{2n+2l}+t_{2n+1+2k}t_{2n+1+2l}\Big)+o(t^2)\\
\phi_{\bar{k}\,\bar{l}}&=&d_n\Big(t_{2k-1}t_{2n+2l}-t_{2k}t_{2n+1+2l}-t_{2l-1}t_{2n+2k}+t_{2l}t_{2n+1+2k}\Big)+o(t^2),\nonumber\\
\tilde{\phi}_{1\,k}&=&d_n\Big(t_{2k}-t_2t_{2k-1}-t_{2n+3}t_{2n+2k}\Big)+o(t^2),\nonumber\\
\tilde{\phi}_{1\,\bar k}&=&d_n\Big(t_{2n+1+2k}+t_2t_{2n+2k}-t_{2k-1}t_{2n+3}\Big)+o(t^2),\nonumber\\
\tilde{\phi}_{k\,l}&=&d_n\Big(t_{2k-1}t_{2l}-t_{2k}t_{2l-1}+t_{2n+2k}t_{2n++1+2l}-t_{2n+1+2k}t_{2n+2l}\Big)+o(t^2),\nonumber\\
\tilde{\phi}_{\bar{k}\,\bar{l}}&=&d_n\Big(t_{2k-1}t_{2n+1+2l}-t_{2k}t_{2n+2l}-t_{2l-1}t_{2n+1+2k}-t_{2l}t_{2n+2k}\Big)+o(t^2),\nonumber
\end{eqnarray}
where $2\leq\alpha\leq n$, $2\leq k\leq n+1$ and $2\leq k<l\leq n+1$.
               Noticing that $F'_*(\epsilon_A)=\frac{\partial F'}{\partial t_A}\Big|_{t=0}$,
               by (\ref{3.6}), (2.1) and (2.2), we get
\begin{eqnarray*}
F'_*(\epsilon_1)&=&\Big(0,\;\tilde{a}_n\,\mbox{E}_1,\;0,\;0,\;0,\;0,\;0\Big),\hspace{0.3cm}\\
F'_*(\epsilon_2)&=&\Big(\tilde{a}_n\,\mbox{E}_2,\;0,\;0,\;0,\;0,\;0,\;0\Big){,}\\
F'_*(\epsilon_{2n+3})&=&\Big(0,\;\tilde{a}_n\,\mbox{E}_2,\;0,\;0,\;0,\;0,\;0\Big),\hspace{0.3cm}\\
F'_*(\epsilon_{2k-1})&=&\Big(\tilde{a}_n\mbox{E}_{2k-1},\;0,\;0,\;\tilde{b}_nd_n\,\mbox{E}_{1\;k},\;0,\;0,\;0\Big),\\
F'_*(\epsilon_{2k})&=&\Big(\tilde{a}_n\mbox{E}_{2k},\;0,\;0,\;0,\;0,\;\tilde{b}_nd_n\,\mbox{E}_{1\;k},\;0\Big),\\
F'_*(\epsilon_{2n+2k})&=&\Big(0,\;\tilde{a}_n\mbox{E}_{2k-1},\;0,\;0,\;\tilde{b}_nd_n\,\mbox{E}_{1\;k},\;0,\;0\Big),\\
F'_*(\epsilon_{2n+1+2k})&=&\Big(0,\;\tilde{a}_n\mbox{E}_{2k},\;0,\;0,\;0,\;0,\;\tilde{b}_nd_n\,\mbox{E}_{1\;k}\Big),
\end{eqnarray*}
where $2\leq k\leq n+1$.
                   Further, we obtain an orthonormal basis of $F'_*(T_pS^{4n+3})$:
\begin{equation*}
e_{A}=\left\{
\begin{array}{lll}
\frac{1}{\tilde{a}_n}F'_*(\epsilon_A),\hspace{0.2cm}A=1,\;2,\;2n+3,\\
{}\\
\frac{1}{\sqrt{\tilde{a}_n^2+\tilde{b}_n^2d_n^2}}F'_*(\epsilon _A),\hspace{0.2cm}A\neq1,\;2,\;2n+3.
\end{array}
\right.
\end{equation*}
                  and an orthonormal basis for the normal space of $F'_*(T_pS^{4n+3})$ in $T_PS^{2n^2+7n+3}$:
\begin{eqnarray*}
e_{4n+4}&=&\Big(-\tilde{b}_n\,\mbox{E}_1,\;0,\;\tilde{a}_n\,\mbox{E}_1,\;0,\;0,\;0\Big),\\
e_{4n+2k+1}&=&\frac{1}{\sqrt{\tilde{a}_n^2+\tilde{b}_n^2d_n^2}}\Big(-\tilde{b}_nd_n\,\mbox{E}_{2k-1},\;0,\;0,
\;\tilde{a}_n\,\mbox{E}_{1\,k},\;0,\;0,\;0\Big),\\
e_{4n+2k+2}&=&\frac{1}{\sqrt{\tilde{a}_n^2+\tilde{b}_n^2d_n^2}}\Big(-\tilde{b}_nd_n\,\mbox{E}_{2k},\;0,\;0,
\;0,\;0,\;\tilde{a}_n\,\mbox{E}_{1\,k},\;0\Big),\\
e_{6n+2k+1}&=&\frac{1}{\sqrt{\tilde{a}_n^2+\tilde{b}_n^2d_n^2}}\Big(\;0,-\tilde{b}_nd_n\,\mbox{E}_{2k-1},\;0,
\;0,\;\tilde{a}_n\,\mbox{E}_{1\,k},\;0,\;0\Big),\\
e_{6n+2k+2}&=&\frac{1}{\sqrt{\tilde{a}_n^2+\tilde{b}_n^2d_n^2}}\Big(\;0,-\tilde{b}_nd_n\,\mbox{E}_{2k},\;0,
\;0,\;0,\;0,\;\tilde{a}_n\,\mbox{E}_{1\,k}\Big),\\
e_{8n+3+\alpha}&=&\Big(0,\hspace{0.2cm}0,\hspace{0.2cm}\mbox{E}_\alpha,
\hspace{0.2cm}0,\hspace{0.2cm}0,\hspace{0.2cm}0,\hspace{0.2cm}0\Big),
\end{eqnarray*}
and
\begin{eqnarray*}
e_{k\,l}=\Big(0,\hspace{0.2cm}0,\hspace{0.2cm}0,\hspace{0.2cm}\mbox{E}_{k\,l},\hspace{0.2cm}0,\hspace{0.2cm}0,\hspace{0.2cm}0\Big),\hspace{0.5cm}e_{\bar k\,\bar l}=\Big(0,\hspace{0.2cm}0,\hspace{0.2cm}0,\hspace{0.2cm}0,\hspace{0.2cm}\mbox{E}_{k\,l},\hspace{0.2cm}0,\hspace{0.2cm}0\Big),\\
\tilde{e}_{k\,l}=\Big(0,\hspace{0.2cm}0,\hspace{0.2cm}0,\hspace{0.2cm}0,\hspace{0.2cm}0,\hspace{0.2cm}\mbox{E}_{k\,l},\hspace{0.2cm}0\Big),
\hspace{0.5cm}\tilde{e}_{\bar{k}\,\bar{l}}=\Big(0,\hspace{0.2cm}0,\hspace{0.2cm}0,\hspace{0.2cm}0,\hspace{0.2cm}0,\hspace{0.2cm}0,\hspace{0.2cm}\mbox{E}_{k\,l}\Big),
\end{eqnarray*}
where $2\leq k\leq n+1$, $2\leq\alpha\leq n$ and $2\leq k<l\leq n+1$.

\vspace{0.3cm}

We set $F'_{AB}=\frac{\partial^2F'}{\partial t_A\partial t_B}\Big|_{t=0}$ and define
\begin{equation*}
H_{AB}=\left\{
\begin{array}{lllll}
\frac{1}{\tilde{a}_n^2}F'_{A\,B},\hspace{0.2cm}A=1,\;2,\;2n+3,\;B=1,\;2,\;2n+3{,}\\
{}\\
\frac{1}{\tilde{a}_n\;\sqrt{a_n^2+\tilde{b}_n^2d_n^2}}F'_{A\,B},\hspace{0.2cm}A=1,\;2,\;2n+3,\;B\neq1,\;2,\;2n+3,\\
{}\\
\frac{1}{\tilde{a}_n^2+\tilde{b}_n^2d_n^2}F'_{A\,B},\hspace{0.2cm}A, B\neq1,\;2,\;2n+3.
\end{array}
\right.
\end{equation*}
                  Then, at $P$, the second fundamental form of $F'$
                  is given in terms of the frame
                  $\{e_A,\;e_\tau,\;e_{k\,l},$ $e_{\bar k\,\bar l},$ $\tilde{e}_{k\,l},\;
                  \tilde{e}_{\bar k\,\bar l}\;|\;1\leq A\leq 4n+3,$ $4n+4\leq\tau\leq 9n+3,$ $2\leq k<l\leq n+1\}$
                   by
\begin{equation*}
h^\tau_{AB}=\big\langle H_{AB},e_\tau\big\rangle,\hspace{0.5cm}h^{(k,\,l)}_{AB}=\big\langle H_{AB},e_{k\,l}\big\rangle,\hspace{0.5cm} h^{(\bar k,\,\bar l)}_{AB}=\big\langle H_{AB},e_{\bar{k}\,\bar{l}}\big\rangle,
\end{equation*}
\begin{equation*}
h^{[k,\,l]}_{AB}=\big\langle H_{AB},\tilde{e}_{k\,l}\big\rangle,\hspace{3.9cm} h^{[\bar k,\,\bar l]}_{AB}=\big\langle H_{AB},\tilde{e}_{\bar{k}\,\bar{l}}\big\rangle.
\end{equation*}
{\ }\\

In summary, we have:

\begin{prop}\label{proposition}
The second fundamental form of $F'$ at the base point $P$ w.r.t. the frame $\{e_A,\;e_\tau,$ $e_{k\,l},e_{\bar k\,\bar l},\;\tilde{e}_{k\,l},$ $\tilde{e}_{\bar k\,\bar l}\;|\;1\leq A\leq 4n+3,4n+4\leq\tau\leq 9n+3,\;2\leq k<l\leq n+1\}$ is given by

\emph{(1)} \hspace{0.3cm} $h^{4n\!+\!4}_{A\,A}=\frac{\tilde{b}_n}{\tilde{a}_n}$,\hspace{0.3cm}   $h^{4n\!+\!4}_{BB}=-\frac{(n+2)\tilde{a}_n\tilde{b}_n}{n(\tilde{a}_n^2+\tilde{b}_n^2d_n^2)}$,\hspace{0.3cm}  $A=1,\;2,\;2n+3$, $B\neq1,\;2,\;2n+3;$

\emph{(2)} \hspace{0.3cm} $h^{4n\!+\!1\!+\!2k}_{1\;2n+2k}=h^{4n\!+\!1\!+\!2k}_{2\;2k}=h^{4n\!+\!1\!+\!2k}_{2n+3\;2n+1+2k}
=\frac{\tilde{b}_nd_n}{\tilde{a}_n^2+\tilde{b}_n^2d_n^2}$,\hspace{0.3cm} $2\leq k\le n+1;$

\emph{(3)} \hspace{0.3cm}$h^{4n\!+\!2\!+\!2k}_{1\;2n+1+2k}=-h^{4n\!+\!2\!+\!2k}_{2\;2k-1}=-h^{4n\!+\! {2}\!+\!2k}_{2n+3\;2n+2k}
=\frac{\tilde{b}_nd_n}{\tilde{a}_n^2+\tilde{b}_n^2d_n^2}$,\hspace{0.3cm} $2\leq k\le n+1;$

\emph{(4)} \hspace{0.3cm}$h^{6n\!+\!1\!+\!2k}_{1\;2k-1}=h^{6n\!+\!1\!+\!2k}_{2\;2n+1+2k}=-h^{6n\!+\!1\!+\!2k}_{2n+3\;2k}
=-\frac{\tilde{b}_nd_n}{\tilde{a}_n^2+\tilde{b}_n^2d_n^2}$,\hspace{0.3cm} $2\leq k\le n+1;$

\emph{(5)} \hspace{0.3cm}$h^{6n\!+\!2\!+\!2k}_{1\;2k}=-h^{6n\!+\!2\!+\!2k}_{2\;2n+2k}=h^{6n\!+\!1\!+\!2k}_{2n+3\;2k-1}
=-\frac{\tilde{b}_nd_n}{\tilde{a}_n^2+\tilde{b}_n^2d_n^2}$,\hspace{0.3cm} $2\leq k\le n+1;$

\emph{(6)} \hspace{0.3cm}$h^{8n\!+\!3\!+\!\alpha}_{2k-1\;2k-1}=h^{8n\!+\!3\!+\!\alpha}_{2k\;2k}=h^{8n\!+\!3\!+\!\alpha}_{2n+2k\;2n+2k}
=h^{8n\!+\!3\!+\!\alpha}_{2n+1+2k\;2n+1+2k}=-\frac{2\tilde{b}_nc_{n,\alpha}}{(n+1-\alpha)(\tilde{a}_n^2+\tilde{b}_n^2d_n^2)}$,\hspace{0.3cm}\\$2\leq \alpha\leq n$, $\alpha< k\le n+1;$

\emph{(7)} \hspace{0.3cm}$h^{8n\!+\!3\!+\!\alpha}_{2\alpha-1\;2\alpha-1}=h^{8n\!+\!3\!+\!\alpha}_{2\alpha\;2\alpha}=h^{8n\!+\!3\!+\!\alpha}_{2n+2\alpha\;2n+2\alpha}
=h^{8n\!+\!3\!+\!\alpha}_{2n+1+2\alpha\;2n+1+2\alpha}=\frac{2\tilde{b}_nc_{n,\alpha}}{\tilde{a}_n^2+\tilde{b}_n^2d_n^2}$,\hspace{0.3cm}$2\leq \alpha\leq n;$


\emph{(8)}\hspace{0.3cm}$h^{(k,l)}_{2k-1\;2l-1}=h^{(k,l)}_{2k\;2l}=h^{(k,l)}_{2n+2k\;2n+2l}=h^{(k,l)}_{2n+1+2k\;2n+1+2l}
=h^{(\bar{k},\bar{l})}_{2k-1\;2n+2l}=-h^{(\bar{k},\bar{l})}_{2k\;2n+1+2l}=\\-h^{(\bar{k},\bar{l})}_{2l-1\;2n+2k}=h^{(\bar{k},\bar{l})}_{2l\;2n+1+2k}=
\frac{\tilde{b}_nd_n}{\tilde{a}_n^2+\tilde{b}_n^2d_n^2},$
\hspace{0.3cm} $2\leq k<l\le n+1;$

\emph{(9)}\hspace{0.3cm}$h^{[k,l]}_{2k-1\;2l}=-h^{[k,l]}_{2k\;2l-1}=h^{[k,l]}_{2n+2k\;2n+1+2l}=-h^{[k,l]}_{2n+1+2k\;2n+2l}
=h^{[\bar{k},\bar{l}]}_{2k-1\;2n+1+2l}=h^{[\bar{k},\bar{l}]}_{2k\;2n+2l}
=\\-h^{[\bar{k},\bar{l}]}_{2l-1\;2n+1+2k}=-h^{[\bar{k},\bar{l}]}_{2l\;2n+2k}=
\frac{\tilde{b}_nd_n}{\tilde{a}_n^2+\tilde{b}_n^2d_n^2},$
\hspace{0.3cm} $2\leq k<l\le n+1;$

\noindent
                       with the same value in the symmetric slot and zero for others.
\end{prop}

\emph{Proof}. By direct computation.
                       Moreover, one can see that $F'$ is minimal.
                       \hfill$\Box$
{\ }\\

\subsection{Type-III Lawson-Osserman sphere}\label{s3.3}
                Let $p_1,p_2\in \textbf{O}$
                written as $p_1=(z_1+z_2\texttt{j})+(z_3+z_4\texttt{j})\texttt{e}$ and
                $p_2=(z_5+z_6\texttt{j})+(z_7+z_8\texttt{j})\texttt{e}$.
                Then,
\begin{eqnarray*}
|p_1|^2-|p_2|^2=\sum\limits_{k=1}^4|z_k|^2-\sum\limits_{k=5}^8|z_k|^2,\\
2\bar{p}_1p_2=\big(a_1+a_2\texttt{j}\big)+\big(a_3+a_4\texttt{j}\big)\texttt{e},
\end{eqnarray*}
where
\begin{eqnarray*}
a_1&=&\bar z_1z_5+z_2\bar z_6+z_3\bar z_7+\bar z_4z_8,\hspace{0.5cm}a_2=\bar z_1z_6-z_2\bar z_5-\bar z_3 z_8+ z_4\bar z_7,\\
a_3&=&\bar z_1z_7+\bar z_2 z_8-z_3\bar z_5- z_4\bar z_6,\hspace{0.5cm}a_4= z_1z_8-z_2z_7+z_3 z_6- z_4 z_5.
\end{eqnarray*}
Set
\begin{equation*}
f_1=|p_1|^2-|p_2|^2,\hspace{0.5cm}f_{1+k}=2\mbox{Re}(a_k),\hspace{0.5cm}f_{5+k}=2\mbox{Im}(a_k),\hspace{0.5cm}1\leq k\leq4,
\end{equation*}
                 Then the third Hopf map $\eta'':S^{15}\longrightarrow S^8$ is given by
                 $z\mapsto\big(f_1,\ldots,f_9\big)$, where $p=(p_1,p_2)$ is identified with $z=(z_1,\ldots,z_8)$.
                 The Type-III Lawson-Osserman sphere is represented by
                 \begin{equation*}
                 F'':S^{15}(1)\longrightarrow S^{24}(1),\hspace{0.5cm}z\mapsto\Big(\sqrt{\frac{28}{45}}\big(\mbox{Re}(z),\mbox{Im}(z)\big),\sqrt{\frac{17}{45}}\big(f_1,f_2\ldots,f_9\big)\Big).
\end{equation*}
                 It is known that $F''$ is also homogeneous.

\vspace{0.3cm}
{For} the second fundamental form of $F''$ at point $P=\Big(\sqrt{\frac{28}{45}}\mbox{E}_1,\;0,\;\sqrt{\frac{17}{45}}\mbox{E}_1\Big)${,}
{we substitute} (2.1) and (2.2) with $m=7$ into $f_k$
{and} obtain
\begin{eqnarray*}\label{3.5}
f_1&=&1-2\sum\limits_{k=5}^8\big(t_k^2+t_{7+k}^2\big)+o(t^2),\\
f_2&=&2\big(t_5+t_2t_6+t_3t_7+t_4t_8+t_9t_{13}+t_{10}t_{14}+t_{11}t_{15}\big)+o(t^2),\\
f_3&=&2\big(t_6-t_2t_5-t_3t_8+t_4t_7-t_9t_{12}-t_{10}t_{15}+t_{11}t_{14}\big)+o(t^2),\\
f_4&=&2\big(t_7+t_2t_8-t_3t_5-t_4t_6+t_9t_{15}-t_{10}t_{12}-t_{11}t_{13}\big)+o(t^2),\\
f_5&=&2\big(t_8-2t_1t_{15}-t_2t_7+t_3t_6-t_4t_5+t_9t_{14}-t_{10}t_{13}+t_{11}t_{12}\big)+o(t^2),\\
f_6&=&2\big(t_{12}-t_2t_{13}-t_3t_{14}+t_4t_{15}+t_6t_{9}+t_{7}t_{10}-t_{8}t_{11}\big)+o(t^2),\\
f_7&=&2\big(t_{13}+t_2t_{12}-t_3t_{15}-t_4t_{14}-t_5t_{9}+t_{7}t_{11}+t_{8}t_{10}\big)+o(t^2),\\
f_8&=&2\big(t_{14}+t_2t_{15}+t_3t_{12}+t_4t_{13}-t_5t_{10}-t_{6}t_{11}-t_{8}t_{9}\big)+o(t^2),\\
f_9&=&2\big(t_{15}+2t_1t_8-t_2t_{14}+t_3t_{13}-t_4t_{12}-t_5t_{11}+t_{6}t_{10}-t_{7}t_{9}\big)+o(t^2).
\end{eqnarray*}
Taking the partial derivative w.r.t. $t_A$, at $t=0$, we have
\begin{eqnarray*}
F''_*(\epsilon_1)=\Big(0,\;\sqrt{\frac{28}{45}}\mbox{E}_1,\;0\Big),\hspace{0.3cm}
F''_*(\epsilon_k)=\Big(\sqrt{\frac{28}{45}}\mbox{E}_k,\;0,\;0\Big),\hspace{0.3cm}
F''_*(\epsilon_{l})=\Big(\sqrt{\frac{28}{45}}\mbox{E}_l,\;0,\;\sqrt{\frac{17}{45}}\mbox{E}_{l-3}\Big),\\
F''_*(\epsilon_{7+k})=\Big(0,\;\sqrt{\frac{28}{45}}\mbox{E}_k,\;0\Big),\hspace{4.2cm}
F''_*(\epsilon_{7+l})=\Big(\;0,\sqrt{\frac{28}{45}}\mbox{E}_l,\;\sqrt{\frac{17}{45}}\mbox{E}_{l+1}\Big),
\end{eqnarray*}
where $2\leq k\leq 3$ and $5\leq l\leq 8$.
                Further, we gain an orthonormal basis of $F''_*(T_pS^{15})$:
\begin{equation*}
e_{A}=\left\{
\begin{array}{lll}
\sqrt{\frac{45}{28}}F''_*(\epsilon_A),\hspace{0.2cm}A=1,\;2,\;3,\;4,\;9,\;10,\;11,\\
{}\\
\sqrt{\frac{15}{32}}F''_*(\epsilon _A),\hspace{0.2cm}A=5,\;6,\;7,\;8,\;12,\;13,\;14,\;15.
\end{array}
\right.
\end{equation*}
                  and an orthonormal basis for the normal space of $F''_*(T_pS^{15})$ in $T_PS^{24}$:
\begin{eqnarray*}
e_{16}&=&\Big(-\sqrt{\frac{17}{45}}\mbox{E}_1,\;0,\;\sqrt{\frac{28}{45}}\mbox{E}_1\Big),\\
e_{12+k}&=&\Big(-\sqrt{\frac{17}{24}}\mbox{E}_k,\;0,\;\sqrt{\frac{7}{24}}\mbox{E}_{k-3}\Big),\\
e_{16+k}&=&\Big(\;0,-\sqrt{\frac{17}{24}}\mbox{E}_l,\;\sqrt{\frac{7}{24}}\mbox{E}_{l+1}\Big).
\end{eqnarray*}
where $5\leq k\leq 8$.

\vspace{0.3cm}

We set $F''_{AB}=\frac{\partial^2F''}{\partial t_A\partial t_B}\Big|_{t=0}$ and define
\begin{equation*}
H_{AB}=\left\{
\begin{array}{lllll}
\frac{45}{28}F''_{A\,B},\hspace{0.2cm}A,B=1,\;2,\;3,\;4,\;9,\;10,\;11,\\
{}\\
\sqrt{\frac{3}{14}}\cdot\frac{15}{8}\;\;F''_{A\,B},\hspace{0.2cm}A=1,\;2,\;3,\;4,\;9,\;10,\;11,\;B=5,\;6,\;7,\;8,\;12,\;13,\;14,\;15,\\
{}\\
\frac{15}{32}F''_{A\,B},\hspace{0.2cm}A, B=5,\;6,\;7,\;8,\;12,\;13,\;14,\;15.
\end{array}
\right.
\end{equation*}
               Then, at $P$, the second fundamental form of $F''$
               is given in terms of the frame $\{e_A,\;e_\tau,|\;1\leq A\leq 15,16\leq\tau\leq 24\}$ by
               $h^\tau_{AB}=\big\langle H_{AB},e_\tau\big\rangle$.
               \\
{\ }\\

In summary, we have:

\begin{prop}\label{proposition}
           The second fundamental form of $F''$ at the base point $P$
           w.r.t. the frame $\{e_A,\;e_\tau\;|$ $1\leq A\leq 15,\;16\leq\tau\leq 24\}$ is given by

\emph{(1)} \hspace{0.3cm} $h^{16}_{1\;1}=h^{16}_{k\;k}=h^{16}_{7+k\;7+k}=\sqrt{\frac{17}{28}}$, $h^{16}_{l\;l}=h^{16}_{l+7\;l+7}=-\frac{\sqrt{119}}{16}$,
$2\leq k\leq  {4}$, $5\leq l\leq 8;$

\emph{(2)} \hspace{0.3cm} $h^{17}_{1\;12}=h^{17}_{2\;6}=h^{17}_{3\;7}=h^{17}_{4\;8}=h^{17}_{9\;13}=h^{17}_{10\;14}=h^{17}_{11\;15}=\frac{\sqrt{85}}{16};$

\emph{(3)} \hspace{0.3cm}$h^{18}_{1\;13}=-h^{18}_{2\;5}=-h^{18}_{3\;8}=h^{18}_{4\;7}=-h^{18}_{9\;12}=-h^{18}_{10\;15}=h^{18}_{11\;14}=\frac{\sqrt{85}}{16};$

\emph{(4)} \hspace{0.3cm}$h^{19}_{1\;14}=h^{19}_{2\;8}=-h^{19}_{3\;5}=-h^{19}_{4\;6}=h^{19}_{9\;15}=-h^{19}_{10\;12}=-h^{19}_{11\;13}=\frac{\sqrt{85}}{16};$

\emph{(5)} \hspace{0.3cm}$h^{20}_{1\;15}=h^{20}_{2\;7}=-h^{20}_{3\;6}=h^{20}_{4\;5}=-h^{20}_{9\;14}=h^{20}_{10\;13}=-h^{20}_{11\;12}=-\frac{\sqrt{85}}{16};$

\emph{(6)} \hspace{0.3cm}$h^{21}_{1\;5}=h^{21}_{2\;13}=h^{21}_{3\;14}=-h^{21}_{4\;15}=-h^{21}_{6\;9}=-h^{21}_{7\;10}=h^{21}_{8\;11}=-\frac{\sqrt{85}}{16};$

\emph{(7)} \hspace{0.3cm}$h^{22}_{1\;6}=-h^{22}_{2\;12}=h^{22}_{3\;15}=h^{22}_{4\;14}=h^{22}_{5\;9}=-h^{22}_{7\;11}=-h^{22}_{8\;10}=-\frac{\sqrt{85}}{16};$

\emph{(8)} \hspace{0.3cm}$h^{23}_{1\;7}=-h^{23}_{2\;15}=-h^{23}_{3\;12}=-h^{23}_{4\;13}=h^{23}_{5\;10}=h^{23}_{6\;11}=h^{23}_{8\;9}=-\frac{\sqrt{85}}{16};$

\emph{(9)} \hspace{0.3cm}$h^{24}_{1\;8}=h^{24}_{2\;14}=h^{24}_{3\;13}=-h^{24}_{4\;12}=-h^{24}_{5\;11}=h^{24}_{6\;10}=-h^{24}_{7\;9}=\frac{\sqrt{85}}{16};$

\noindent
                     with the same value in the symmetric slot
                     and zero for others.
\end{prop}

\emph{Proof}.
                      By computation. Moreover, one can see that $F''$ is minimal.
                       \hfill$\Box$
{\ }\\{\ }\\
\section{On the area-minimizing property}\label{sec4}

\subsection{Lawlor's curvature criterion}
For completeness, we briefly recall Lawlor's curvature criterion for proving a minimal cone to be area-minimizing.
For further details readers are referred to \cite{[L91]}.

            Let $\Sigma$ be a smooth $n$-dimensional submanifold of the unit sphere $S^N$ and
                  \[
                  \mathcal{C}\Sigma=\{tx:t\in[0,\infty) \text{ and } x\in \Sigma\}.
                  \]
            Fix $p\in\Sigma$.
            A \emph{normal geodesic} of length $\ell$ is an arc of a great circle $\gamma$
            which is perpendicular to $\Sigma$ at its starting point $\gamma(0)=p$.
            We call $\gamma$ an \emph{open normal geodesic}
            if we leave off the endpoint $\gamma(\ell)$.
            Let $U_p(\ell)$ be the union of points of open normal geodesics from $p$ of length $\ell$.
            Then \emph{normal wedge} $W_p(\ell)$ is defined to be $\mathcal{C}U_p(\ell)-\{0\}$.
            The \emph{normal radius} of $\mathcal{C}\Sigma$ at a point $p\in\Sigma$ is the largest $\ell_p$
            such that
            $W_p(\ell_p)$ intersects $\mathcal{C}\Sigma$ only in the ray $\overrightarrow{op}$.

             Suppose $p\in\Sigma$ and $\nu$ is a unit vector in the normal space $T^\perp_p\Sigma$.
             Let $(r,\theta)$ be the polar coordinate of the plane spanned by $\overrightarrow{op}$ and $\nu$.
             A \emph{projection curve} $\gamma_p$, if exists, satisfies

\begin{equation*}
(\mbox{ODE})\hspace{1cm}\left\{
\begin{array}{lll}
\frac{dr}{d\theta}=r\sqrt{r^{ {2n+2}}\cos^{2n}\theta\inf\limits_{\nu\in T^\perp_p\Sigma,\;|\nu|=1}\Big(\det\big(\mbox{I}-\tan\theta\;(h^\nu_{AB})\big)\Big)^2-1},\\
{}\\
r(0)=1,\\
\end{array}
\right.
\end{equation*}
                where $(h^\nu_{AB})$ is the matrix of the second fundamental form of $\Sigma$ at $p$, in the normal direction $\nu$.
                The existence of the ODE relies on the size of second fundamental form and 
                the dimension of $\mathcal{C}\Sigma$.
           If $\gamma_p$ exists, either $\frac{dr}{d\theta}$ vanishes at some positive $\theta(p)$,
           or $r$ goes to infinity as $\theta$ approaches some finite value $\theta_0(p)$.
           In the latter case, we call the smallest $\theta_0(p)$ the vanishing angle at $p$.
           Let $\Gamma_p$ be the rotated surface generated by $\gamma_p$ in $W_p(\theta_0(p))$.
           Then we define $\Pi_p$ by sending $\Gamma_p$ to $p$ and requiring $\Pi_{p}(tz)=t\Pi_p(z)$ for $t>0$ and $z\in\Gamma_p$.
           If $\{W_p(\theta_0(p)):p\in\Sigma\}$ do not intersect,
           we assemble $\{\Pi_p:p\in\Sigma\}$ together and
           extend it to a global
           \emph{retraction} $\Pi:\textbf{R}^{N+1}\longrightarrow\mathcal{C}\Sigma$
           which equals $\Pi_p$ in $W_p(\theta_0(p))$ and
           collapses everything else to $0$.
           It can be guaranteed by (ODE) that $\Pi$ is a continuously area-noincreasing projection to $\mathcal{C}\Sigma$.

\vspace{0.1cm}
By using the retraction $\Pi$, Lawlor proved
\begin{theo}\label{theorem4.1}
\emph{(\textbf{Lawlor's curvature criterion} \cite{[L91]})}
Let $\Sigma$ be a smooth $n$-dimensional submanifold of unit sphere $S^N$.
            Suppose 
            that the vanishing angle $\theta_0(p)$ exists for every $p\in \Sigma$ and that
            $\ell_0=\min\limits_{p\in\Sigma} \ell_p \geq 2\max\limits_{p\in\Sigma} \theta_0(p)$
            which ensures that $\{W_p(\theta_0(p)):p\in\Sigma\}$ do not intersect.
            Then $\mathcal{C}\Sigma$ is area-minimizing
            (in the sense of mod 2 when $\Sigma$ is nonorientable).
\end{theo}
\emph{Remark}.
                            Lawlor made a table (page 20-21 in \cite{[L91]}) of estimated vanishing angles for $\dim \mathcal{C}\Sigma\leq12$
                            and
                                    $\mathcal S^2$ where
                                    \begin{equation}\label{defCalS}
 {\mathscr S}=\max\limits_{p\in\Sigma}\Big(\sup\limits_{\nu\in T^\perp_p\Sigma,\;|\nu|=1}\big(\sum\limits_{A,B}(h_{AB}^\nu)^2\big)^{\frac{1}{2}}\Big).
\end{equation}
                            He used the control
                                  \begin{equation}\label{nstrctrl}
                                  \inf\limits_{\nu\in T^\perp_p\Sigma,\;|\nu|=1}\Big(\det\big(\mbox{I}-t\;(h^\nu_{AB})\big)\Big)
                                  {>} (1-\mathcal S t) e^{\mathcal S t}
                                  \end{equation}
                              for $\dim \mathcal C \Sigma=12$
                              and a more accurate lower bound $F(\mathcal S, t, \dim \mathcal C \Sigma)$ for $\dim \mathcal C \Sigma<12$.
                                                            By $V(m, {\mathscr S})$ we mean the estimated vanishing angle based on \eqref{nstrctrl}
                            for $m=\dim \mathcal{C}\Sigma\geq12$ and $ {\mathscr S}$.
                             When $m>12$, Lawlor proved the following nice property
                            \begin{equation}\label{4.1}
\tan\Big(V(m,\frac{m}{12} {\mathscr S})\Big)<\frac{12}{m}\tan\Big(V(12, {\mathscr S})\Big).
\end{equation}
                                                       Moreover, we remark that
                            \begin{equation}\label{mono}
                            V(m,a)<V(m,b) \text{ for } a<b.
                            \end{equation}
{\ }\\{\ }

\subsection{Proof of the main theorem}

      Let $F$, $F'$ and $F''$ be the Lawson-Osserman spheres constructed in \S \ref{s3.1}, \S\ref{s3.2} and \S\ref{s3.3} respectively.
       Then, in this subsection, we prove our
       
 \vspace{0.3cm}
     \textbf{Main Theorem}. \emph{The minimal cones $\mathcal{C}F$, $\mathcal{C}F'$ and $\mathcal{C}F''$ are area-minimizing.}
       \vspace{0.3cm}
    {\ }\\   
\emph{Proof}. To reduce redundance, we present a complete proof only for $\mathcal{C}F$.
                  We will show

(1) Any normal line through $P$ intersects $\mathcal{C}F$ only at $P$, i.e., the normal radius $\ell_0\geq\frac{\pi}{2}$;

(2) The vanishing angle $\theta_0<\frac{\pi}{4}$;

\noindent and the theorem follows by Lawlor's criterion.

                 Since $F$, $F'$ and $F''$ are homogeneous,
                  it is sufficient to do calculations at the base point $P$.
                We verify (1) first.
                Let
                {$X-P$ be a} normal vector through $P$.
                Then, according to (\ref{3.3-1}), $X$ can be written as
\begin{eqnarray*}
           X=P&+&\lambda_{2n+2}\,e_{2n+2}+\sum\limits_{k=2}^{n+1}\Big(\lambda_{2n+1+k}\;e_{2n+1+k}+\lambda_{3n+1+k}\;e_{3n+1+k}\Big)\\
&+&\sum\limits_{\alpha=2}^{n}\lambda_{4n+1+\alpha}\;e_{4n+1+\alpha}+\sum\limits_{2\leq k<l\leq n+1}\Big(\lambda_{k\;l}\;e_{k\;l}+\lambda_{\bar k\;\bar l}\;e_{\bar k\;\bar l}\Big).
\end{eqnarray*}
                    In terms of blocks, 
                    we write $X=\big(\xi,\eta,\mu,\varsigma,\tau\big)$ and (\ref{3.3-1}) gives
\begin{eqnarray*}
\xi_1=a_n-b_n\lambda_{2n+2},\hspace{2.2cm}\xi_k=-\frac{b_nd_n}{\sqrt{a_n^2+b_n^2d_n^2}}\lambda_{2n+1+k},\hspace{0.3cm} 2\leq k\leq n+1,\\
\eta_1=0,\hspace{4.3cm} \eta_k=\frac{b_nd_n}{\sqrt{a_n^2+b_n^2d_n^2}}\lambda_{3n+1+k},\hspace{0.3cm} 2\leq k\leq n+1,\\
\mu_1=b_n+a_n\lambda_{2n+2},\hspace{4.8cm}\mu_\alpha=\lambda_{4n+1+\alpha},\hspace{0.3cm}2\leq \alpha\leq n,\\
\varsigma_{1\;l}=\frac{a_n}{\sqrt{a_n^2+b_n^2d_n^2}}\lambda_{2n+1+l},\hspace{0.3cm}2\leq l\leq n+1,\hspace{0.8cm}\varsigma_{k\;l}=\lambda_{k\;l},\hspace{0.3cm} 2\leq k<l\leq n+1,\\
\tau_{1\;l}=\frac{a_n}{\sqrt{a_n^2+b_n^2d_n^2}}\lambda_{3n+1+l},\hspace{0.3cm}2\leq l\leq n+1,\hspace{0.8cm}\tau_{k\;l}=\lambda_{\bar k\;\bar l},\hspace{0.3cm}2\leq k<l\leq n+1.
\end{eqnarray*}
Assume $\lambda_{2n+3}\neq 0$, then
\begin{equation}\label{4.2}
\frac{\varsigma_{1\;2}}{\xi_2}=-\frac{a_n}{b_nd_n}.
\end{equation}
If $X\in\mathcal{C}F$, we have
\begin{eqnarray}\label{4.3}
\varsigma_{1\;2}=\frac{\sqrt{|\xi|^2+|\eta|^2}}{a_n}\,b_nd_n\;\mbox{Re}(z_1\bar z_2)
=\frac{b_nd_n}{a_n\sqrt{|\xi|^2+|\eta|^2}}\,\xi_1\xi_2.
\end{eqnarray}
Combining (\ref{4.2}) and (\ref{4.3}), we obtain
\begin{equation}\label{4.4}
\frac{\xi_1}{\sqrt{|\xi|^2+|\eta|^2}}=-\frac{a_n^2}{b_n^2d_n^2}
\end{equation}
which implies
\begin{equation}\label{4.5}
a_n-b_n\lambda_{2n+2}<0.
\end{equation}
On the other hand, from Lemma \ref{lemma3.1} and (\ref{4.4}), we have
\begin{eqnarray}\label{4.6}
b_n+a_n\lambda_{2n+2}=\mu_1&=&\frac{\sqrt{|\xi|^2+|\eta|^2}}{a_n}\;b_nc_{n,1}\;\Big[\frac{\xi_1^2}{|\xi|^2+|\eta|^2}-\frac{1}{n}\sum\limits_{k=2}^{n+1}
\frac{\xi_k^2+\eta_k^2}{|\xi|^2+|\eta|^2}\Big]\nonumber\\
&=&\frac{\sqrt{|\xi|^2+|\eta|^2}}{a_n}\;b_nc_{n,1}\;\Big[\frac{\xi_1^2}{|\xi|^2+|\eta|^2}-\frac{1}{n}\Big(1-\frac{\xi_1^2}{|\xi|^2+|\eta|^2}\Big)\Big]
\nonumber\\
&=&\frac{\sqrt{|\xi|^2+|\eta|^2}}{a_n}\;b_nc_{n,1}\;\Big[\frac{n+1}{n}\;\frac{1}{(2n+3)^2}-\frac{1}{n}\Big]<0.
\end{eqnarray}
                      So (\ref{4.5}) and (\ref{4.6}) together lead to a contradiction. 
                      Thus $\lambda_{2n+3}=0$.

                       Similarly, one can show $\lambda_{2n+1+l}=\lambda_{3n+1+k}=0$ for $3\leq l\leq n+1$ and $2\leq k\leq n+1$.
                       Consequently, $\xi_k=\eta_k=0$ for $2\leq k\leq n+1$.
                       Note that $X,P\in\mathcal CF$ and their first $(2n+2)$ components form parallel vectors.
                     By the geometric structure of $\mathcal CF$, this implies that $\overrightarrow{OX} \sslash \overrightarrow{OP}$.
                       Since $\overrightarrow{OX} = \overrightarrow{OP}+ \overrightarrow{PX}$ with
                       $\overrightarrow{OP}\perp \overrightarrow{PX}$,
                       it follows that
                       $\overrightarrow{PX}=\overrightarrow{0}$,
                       i.e., $X=P$.
                       Hence the normal radius of $\mathcal{C}F$ is  {pointwise at least} $\frac{\pi}{2}$.

        {\ }

                   In order to have estimates (2) on vanishing angles,
                   we need to figure out $\mathcal S$ for our cases.
                   {Suppose}
\begin{eqnarray*}
\nu_0&=&\lambda_{2n+2}\,e_{2n+2}+\sum\limits_{k=2}^{n+1}\Big(\lambda_{2n+1+k}\;e_{2n+1+k}+\lambda_{3n+1+k}\;e_{3n+1+k}\Big)\\
&{}&+\sum\limits_{\alpha=2}^{n}\lambda_{4n+1+\alpha}\;e_{4n+1+\alpha}+\sum\limits_{2\leq k<l\leq n+1}\Big(\lambda_{k\;l}\;e_{k\;l}+\lambda_{\bar k\;\bar l}\;e_{\bar k\;\bar l}\Big)
\end{eqnarray*}
is a unit normal vector,
       {
      such that $$\|\left(h^{\nu_0}_{AB}\right)\|=\mathcal S=\max_{\mu\in T_P^\perp \mathcal{C}F,\ |\mu|=1}\|\left(h^\mu_{AB}\right)\|.$$
      }

            Let $I$ be the set of indices of normal basis \eqref{3.3-1}.
            According to behaviors of second fundamental form in normal directions,
            we split $I$ into two parts:
                  \[
                 \mathscr{A}=\{2n+2\}\bigcup \{4n+1+\alpha\; |\ 2\leq \alpha \leq n\}\ \ \    and \ \ \
                  \mathscr{B}=I - \mathscr{A}.
                   \]
           {Then, by Proposition \ref{hforF},

           $(\star 1)$ If $\tau\in\mathscr{A}$, then $(h_{AB}^\tau)$ is purely diagonal, and

           $(\star 2)$ If $\tau\in\mathscr{B}$, then  $\{h_{AA}^\tau\}$ are all zero.
           Moreover, for different $\tau,\mu\in \mathscr{B}$,
           $h^\tau_{AB}\cdot h^\mu_{AB}=0,\ \forall A, B$.
      }

{\ }

             By direct computations, we have
                   \[
                   \|\left(h_{AB}^{2n+2}\right)\|^2=\frac{n(2n+3)}{2(n+1)}+\frac{2n+3}{4(n+1)}  = \frac{(2n+1)(2n+3)}{4(n+1)}
                  \]
                  \[
                    \|\big(h_{AB}^{2n+1+k}\big)\|^2=\|\big(h_{AB}^{3n+1+k}\big)\|^2= \frac{(2n+1)(2n+3)}{4(n+1)(n+2)}\ \ \ for \ \ \ 2\leq k\leq n+1
                   \]
                    \[
                    \|\left(h_{AB}^{4n+1+\alpha}\right)\|^2=\frac{(2n+1)(2n+3)}{2(n+1)(n+2)}\ \ \ for \ \ \ 2\leq \alpha\leq n
                   \]
                    \[
                   \|\big(h_{AB}^{(k,l)}\big)\|^2=\|\big(h_{AB}^{(\overline{k},\overline{l})}\big)\|^2=\frac{(2n+1)(2n+3)}{2(n+1)(n+2)} \ \ \ for \ \ \ 2\leq k< l\leq n+1
                  \]
{\ }\\
             We shall show that $\nu_0=\pm e_{2n+2}$ and ${\mathscr S^2}={\frac{(2n+1)(2n+3)}{4(n+1)}}$.
             The procedure consists of two steps.
\\{\ }

             Step 1: To show $\lambda_\tau=0$ for $\tau\in \mathscr B$.
                         Write $\nu_0=c\cdot \nu_0^\mathscr{A}+s\cdot \nu_0^\mathscr{B}$
                          where $\nu_0^\mathscr{A}, \nu_0^\mathscr{B}$, if not zero, are the normalized unit vectors
                          of projections of $\nu_0$ in span$\{e_\tau{\;|\;\tau\in \mathscr A}\}$ and span$\{e_\tau{\;|\;\tau\in \mathscr B}\}$ respectively,
                          and where $c,s$ stand for $\cos t$ and $\sin t$ for some real number $t$.
                          (The case of either projection of $\nu_0$ being zero is easy to handle.)
                          Then
                          \[
                           \|\big(h_{AB}^{\nu_0}\big)\|^2=
                           \|\big(c\cdot h_{AB}^{\nu_0^\mathscr A} + s\cdot h_{AB}^{\nu_0^\mathscr B}\big)\|^2
                           = c^2\cdot \|\big(h_{AB}^{\nu_0^\mathscr A}\big)\|^2 + s^2\cdot \|\big(h_{AB}^{\nu_0^\mathscr B}\big)\|^2
                          \]
                         By $(\star 2)$, it follows $\|\big(h_{AB}^{\nu_0^\mathscr B}\big)\|^2\leq \frac{(2n+1)(2n+3)}{2(n+1)(n+2)}$.
                         Now, if $\|\big(h_{AB}^{\nu_0^\mathscr A}\big)\|^2$ were
                         bounded from the above by the same number,
                         then
                            \[
                            \|\big(h_{AB}^{\nu_0}\big)\|^2\leq \frac{(2n+1)(2n+3)}{2(n+1)(n+2)}<\|\big(h_{AB}^{2n+2}\big)\|^2
                            \]
                         contradicting with our choice of $\nu_0$.
                         Hence,
                              \[
                               \|\big(h_{AB}^{\nu_0^{\mathscr B}}\big)\|^2\leq \frac{(2n+1)(2n+3)}{2(n+1)(n+2)}<\|\big(h_{AB}^{\nu_0^{\mathscr A}}\big)\|^2
                               \]
                     and consequently
                               \[
                           \|\big(h_{AB}^{\nu_0}\big)\|^2
                           \leq
                           \|\big(h_{AB}^{\nu_0^\mathscr A}\big)\|^2
                           \Rightarrow
                            \|\big(h_{AB}^{\nu_0}\big)\|^2
                           =
                           \|\big(h_{AB}^{\nu_0^\mathscr A}\big)\|^2
                           \Rightarrow
                           \nu_0=\nu_0^\mathscr{A},\ i.e.\ \lambda_\tau=0\ \ \ for\ \ \ \tau\in \mathscr B.
                          \]
                      {\ }

               Step 2: Further to prove $\lambda_\tau=0$ for $\tau\in \mathscr{A}- \{2n+2\}$.
                       We shall deduce the statement by induction. 
                      Write $\nu_0=c\cdot E_{5n+1}+s\cdot e_{5n+1}$
                      where $s=\sin t=\lambda_{5n+1}$ for some $|t|\leq \frac{\pi}{2}$,
                      and where $E_{5n+1}$ is the normalized unit vector of $\nu_0- \lambda_{5n+1}\cdot e_{5n+1}$
                      (nonzero, otherwise contradicting with the choice of $\nu_0$).
                      Note that the nonzero elements of $\left(h_{AB}^{5n+1}\right)$ are
                        \[
                        -h_{n\; n}^{5n+1}
                        =h_{n+1\; n+1}^{5n+1}
                        =-h_{2n\; 2n}^{5n+1}
                        =h_{2n+1\; 2n+1}^{5n+1}
                        =-\frac{2b_nc_{n,n}}{a_n^2+b_n^2d_n^2},
                        \]
                      and meanwhile, by $(\star 1)$ and Proposition \ref{hforF}, that
                      all nonzero elements of $\big(h_{AB}^{E_{5n+1}}\big)$ distribute in its diagonal and
                            \begin{equation*}\label{same}
                        h_{n\; n}^{E_{5n+1}}
                        =h_{n+1\; n+1}^{E_{5n+1}}
                        =h_{2n\; 2n}^{E_{5n+1}}
                        =h_{2n+1\; 2n+1}^{E_{5n+1}}.
                        \end{equation*}
                      These nice distributions support
                      \begin{equation}\label{split}
                           \|\big(h_{AB}^{\nu_0}\big)\|^2
                           = \|\big(c\cdot h_{AB}^{E_{5n+1}} + s\cdot h_{AB}^{5n+1}\big)\|^2
                           = c^2\cdot \|\big(h_{AB}^{E_{5n+1}}\big)\|^2 + s^2\cdot \|\big(h_{AB}^{5n+1}\big)\|^2.
                      \end{equation}
                      By the same argument in Step 1, it follows that $s=\lambda_{5n+1}=0$.

                      Assume that $\lambda_{5n+1}=\lambda_{5n}=\cdots=\lambda_{4n+r+1}=0$ for some $3\leq r\leq n$.
                      We aim to have $\lambda_{4n+r}=0$.
                      Similarly,
                             write $\nu_0=c\cdot E_{4n+r}+s\cdot e_{4n+r}$
                      where $s=\sin t=\lambda_{4n+r}$ for some $|t|\leq \frac{\pi}{2}$,
                      and where $E_{4n+r}=k_{1}e_{2n+2}+\sum\limits_{i=2}^{r-2}k_i e_{4n+i+1}$ is the normalized unit vector of $\nu_0- \lambda_{4n+r}\cdot e_{4n+r}$.
                      Observe, from Proposition \ref{hforF}, that
                       \begin{equation}
                   \left( h_{A\; B}^{{2n+2}}\right)                  =
                                   \diag\left(\frac{b_n}{a_n}, -\frac{(n+2)a_nb_n}{n(a_n^2+b_n^2d_n^2)}, \cdots, -\frac{(n+2)a_nb_n}{n(a_n^2+b_n^2d_n^2)}
                                  \right),
                        \end{equation}
                      and that, for $2\leq \alpha\leq n$,
                      \begin{equation}
                      \begin{split}
                   \left( h_{A\; B}^{{4n+\alpha+1}}\right) &\ \ \ \ \ \ \ \ \ \ \ \ \ \ \ \ \  \ \ {\small \underset{\shortdownarrow}{{\mathbf{\alpha}}\textbf{-th}}}\\
                 =
                                   \diag\Bigg(0, \
                                    & 0, \cdots, 0, \frac{2b_nc_{n,\alpha}}{a_n^2+b_n^2d_n^2}, -\frac{2b_nc_{n,\alpha}}{(n\!+\!1-\alpha)(a_n^2+b_n^2d_n^2)}, \cdots, -\frac{2b_nc_{n,\alpha}}{(n\!+\!1-\alpha)(a_n^2+b_n^2d_n^2)},\\
                                     & 0, \cdots, 0, \frac{2b_nc_{n,\alpha}}{a_n^2+b_n^2d_n^2}, -\frac{2b_nc_{n,\alpha}}{(n\!+\!1-\alpha)(a_n^2+b_n^2d_n^2)} \cdots, -\frac{2b_nc_{n,\alpha}}{(n\!+\!1-\alpha)(a_n^2+b_n^2d_n^2)}
                                  \Bigg),
                                   \end{split}
\end{equation}
                      where the first row includes $n+1$ elements and the second $n$ terms.
                      It follows
                      \begin{equation}
                      \begin{split}
                   \left( h_{A\; B}^{E_{4n+r}}\right) &
                   \ \  \ \ \  \ \ \ \ \ \  \  \ \ \  \ \ \ \ \ \  \ \ \ \  \ \ \ \ \ \ \ \ \ \ \ \  \ \ \ \ \  \ \ \ \ \ \ \ \  \ \ \ \ \ \  \ \ \ \ \  \ {\small \underset{\shortdownarrow}{\textbf{(r-1)-th}}}\\
                 =
                                   \diag\Bigg(k_1\frac{b_n}{a_n}, & \;
                                     -\frac{k_1(n+2)a_nb_n}{n(a_n^2+b_n^2d_n^2)}+\frac{2k_2b_nc_{n,2}}{a_n^2+b_n^2d_n^2},\;
                                     \cdots,\; [*],\;\; \circledast,\; \circledast,\; \cdots,\; \circledast,\\
                                       &\  -\frac{k_1(n+2)a_nb_n}{n(a_n^2+b_n^2d_n^2)}+\frac{2k_2b_nc_{n,2}}{a_n^2+b_n^2d_n^2},\;
                                       \cdots,\; [*],\;\; \circledast,\; \circledast,\; \cdots,\; \circledast
                                  \Bigg),
                                          \end{split}
\end{equation}
                     where all $\circledast$s represent the same number.
                      Hence,
                       \begin{equation}
                      \begin{split}
                   &\left( h_{A\; B}^{\nu_0}\right) =
                                   \diag\Bigg(\left[k_1\frac{b_n}{a_n}\right]c,
                  \\
                 & \cdots, \left[*\right]c,  \circledast c+\left[\frac{2b_nc_{n,r-1}}{a_n^2+b_n^2d_n^2}\right]s,
                                                                        \circledast c-\left[\frac{2b_nc_{n,r-1}}{(n\!+\!2-r)(a_n^2+b_n^2d_n^2)}\right]s,\cdots,
                                                                         \circledast c-\left[\frac{2b_nc_{n,r-1}}{(n\!+\!2-r)(a_n^2+b_n^2d_n^2)}\right]s,\\
                                               &                \cdots, \left[*\right]c,  \circledast c+\left[\frac{2b_nc_{n,r-1}}{a_n^2+b_n^2d_n^2}\right]s,
                                                                        \circledast c-\left[\frac{2b_nc_{n,r-1}}{(n\!+\!2-r)(a_n^2+b_n^2d_n^2)}\right]s,\cdots,
                                                                         \circledast c-\left[\frac{2b_nc_{n,r-1}}{(n\!+\!2-r)(a_n^2+b_n^2d_n^2)}\right]s\Bigg).
                                          \end{split}
\end{equation}
                           Consequently,
                           \begin{equation}
                      \begin{split}
                   &\|\left( h_{A\; B}^{\nu_0}\right)\|^2 \\
                   =\; &
                   c^2\cdot \|\big(h_{A\; B}^{E_{4n+r}}\big)\|^2+s^2\cdot  \|\big(h_{A\; B}^{{4n+r}}\big)\|^2
                   + 4cs\cdot
                                  \left[
                                  \left(\frac{2b_nc_{n,r-1}}{a_n^2+b_n^2d_n^2}\right)\circledast
                                  -(n+2-r)\frac{2b_nc_{n,r-1}}{(n\!+\!2-r)(a_n^2+b_n^2d_n^2)}\circledast
                                   \right]\\
                     =\; &
                     c^2\cdot \|\big(h_{A\; B}^{E_{4n+r}}\big)\|^2+s^2\cdot  \|\big(h_{A\; B}^{{4n+r}}\big)\|^2.
                                                          \end{split}
\end{equation}
                     Similarly, repeating the argument in Step 1 confirms that $\lambda_{4n+r}=0$.
                     Thus, by induction, $\nu_0$ has to be $\pm e_{2n+2}$.

{\ }

                      Now we figure out ${\mathscr S^2}={\frac{(2n+1)(2n+3)}{4(n+1)}}<n+1=\frac{1}{2}\dim(\mathcal{C}F)$.
                      It can be seen from Lawlor's table that,
                      for $2\leq n\leq 5$,
                      vanishing angle $\theta_0$ exists and is less than $45^\circ$.
                      For $n>5$, namely $m=\dim(\mathcal{C}F)=2n+2>12$,
                      by \eqref{4.1} and \eqref{mono} we have
                      \begin{equation}
                      \begin{split}
                   \tan\Big(V(m,\sqrt{\frac{m}{2}})\Big)
                   &=
\tan\Big(V(m,\frac{m}{12}\sqrt{\frac{72}{m}})\Big)
                <
                \frac{12}{m}\tan\Big(V(12,\sqrt{\frac{72}{m}})\Big)\\
                &<
\tan\Big(V(12,\sqrt{6})\Big)
                <\tan 8.36^\circ <1.
                                   \end{split}
\end{equation}
        Hence, $\theta_0$ exists and is less than $\frac{\pi}{4}$ as well.
        Combined with the result in \cite{[HL82]} (about the classical coassociative Lawson-Osserman cone in $\textbf{R}^7$ when our $n=1$), the proof completes.
        \hfill$\Box$

\vspace{0.3cm}

\emph{Remarks}. (a) We would like to point out that the coassociative cone's being area-minimizing cannot be verified following Lawlor's curvature criterion.
Note that
                                \begin{equation}
                                  \det\left(\mbox{I}-t\;(h^{4}_{AB})\right)
                                  \geq
                                  \inf\limits_{\nu\in T^\perp_p\Sigma,\;|\nu|=1}\Big(\det\big(\mbox{I}-t\;(h^\nu_{AB})\big)\Big)
                                 \geq
                                 F\left(\sqrt{\frac{15}{8}},t,3\right)
                                  \end{equation}
                             where $F(\cdot,\cdot,\cdot)$ is the control which Lawlor used for his table of vanishing angles for cones of dimensions below $12$.
                             In this concrete case equalities are attained for all $t$
                             but the case $(\dim, \mathcal S^2)=(4,\frac{15}{8})$ supports no vanishing angle
                             by careful numerical computation.

{\ }\\{\ }

        (b) Type II enjoys similar properties as Type I.
        From \eqref{3.6}, we get
                     \[
                   \|\big(h_{AB}^{4n\!+\!4}\big)\|^2=\frac{n(4n+5)}{2(n+1)}+\frac{3(4n+5)}{8(n+1)}  = \frac{(4n+3)(4n+5)}{8(n+1)}
                  \]
                  \[
                    \|\big(h_{AB}^{4n+1+2k}\big)\|^2=\|h_{AB}^{4n+2+2k}\|^2= \frac{3(4n+3)(4n+5)}{16(n+1)(n+2)}\ \ \ for \ \ \ 2\leq k\leq n+1
                   \]
                   \[
                    \|\big(h_{AB}^{6n+1+2k}\big)\|^2=\|h_{AB}^{6n+2+2k}\|^2= \frac{(4n+3)(4n+5)}{4(n+1)(n+2)}\ \ \ for \ \ \ 2\leq k\leq n+1
                   \]
                    \[
                    \|\big(h_{AB}^{8n+3+\alpha}\big)\|^2=\frac{(4n+3)(4n+5)}{4(n+1)(n+2)}\ \ \ for \ \ \ 2\leq \alpha\leq n
                   \]
                    \[
                   \|\big(h_{AB}^{(k,l)}\big)\|^2=
                    \|\big(h_{AB}^{(\overline{k},\overline{l})}\big)\|^2=
                     \|\big(h_{AB}^{[k,l]}\big)\|^2=
                     \|\big(h_{AB}^{[\overline{k},\overline{l}]}\big)\|^2
                   =\frac{(4n+3)(4n+5)}{4(n+1)(n+2)} \ \ \ for \ \ \ 2\leq k< l\leq n+1
                  \]
             The same idea shows that $\mathcal S^2=\frac{(4n+3)(4n+5)}{8(n+1)}<\frac{1}{2}\dim \mathcal CF'$
             and therefore that $\mathcal{C}F'$ are area-minimizing for all $n\geq 1$.

{\ }

              (c) A similar calculation for Type III gives $\mathcal S^2=\frac{51}{4}+\frac{119}{32}<17$.
              By Lawlor's table of vanishing angles, \eqref{4.1} and \eqref{mono}, it follows that $V(16, \mathcal S)<V(16, \sqrt {17})<V(12, \sqrt {17})\approx 10.11^\circ<\frac{\pi}{4}$.
              So, $\mathcal{C}F''$ turns out to be area-minimizing as well.

\vspace{1cm}

\textbf{Acknowledgments}.
The authors are grateful to Professors Yuanlong Xin, 
Chiakuei Peng,
Zizhou Tang for their helps and constant encouragements during the preparation of this paper,
and the referee for helpful comments to improve the preliminary version.
The third author is indebted to Professor H. Blaine Lawson, Jr. for drawing our attention to \cite{[L91]}.
This work was supported in part by
NSFC (Grant Nos. 11471299, 11471078, 11622103, 11526048, 11601071),
the Fundamental Research Funds for the Central Universities, and the SRF for ROCS, SEM.


{\ }

%




\begin{bibdiv}
\begin{biblist}



\bib{Almgren}{article}{
    author={Almgren, F.J.}
    title={Some interior regularity theorems for minimal surfaces and an extension of Bernstein's theorem},
    journal={Ann. Math.},
    volume={84},
    date={1966},
    pages={277--292},
}



\bib{BDG}{article}{
    author={Bombieri, E.},
     author={De Giorgi, E.},
      author={Giusti, E.},
    title={Minimal cones and the Bernstein problem},
    journal={Invent. Math.},
    volume={7},
    date={1969},
    pages={243--268},
}


\bib{[DW71]}{article}{
    author={do Carmo, M.P.},
       author={Wallach, N.R.},
    title={Minimal immersions of spheres into spheres},
    journal={Ann. Math.},
    volume={93},
    date={1971},
    pages={43--62},
}



\bib{Cheng}{article}{
    author={Cheng, B.N.},
    title={Area-minimizing cone-type surfaces and coflat calibrations},
    journal={Indiana Univ. Math. J.},
    volume={37},
    date={1988},
    pages={505--535},
}

\bib{F}{book}{
    author={Federer, H.},
    title={Geometric Measure Theory},
    place={Springer-Verlag, New York},
    date={1969},
}




\bib{DeG}{article}{
    author={De Giorgi, E.},
    title={Una estensione del teorema di Bernstein},
    journal={Ann. Scuola Norm. Sup. Pisa},
    volume={19},
    date={1965},
    pages={79--85},
}


%
\bib{FK}{article}{
    author={Ferus, D.}
    author={Karcher, H.}
    title={Non-rotaional minimal spheres and minimizing cones},
    journal={Comment. Math. Helv.},
    volume={60},
    date={1985},
    pages={247--269},
}


\bib{Fleming}{article}{
    author={Fleming, W.H.},
    title={On the oriented Plateau problem},
    journal={Rend. Circ. Mat. Palermo},
    volume={11},
    date={1962},
    pages={69--90},
}




\bib{HS}{article}{
    author={Hardt, R.},
    author={Simon, L.},
    title={Area-minimizing hypersurfaces with isolated singularities},
    journal={J. Reine. Angew. Math.},
    volume={362},
    date={1985},
    pages={102--129},
}



\bib{[HL82]}{article}{    author={Harvey, F.R.},    author={{Lawson, Jr.}, H.B.},    title={Calibrated geometries},    journal={Acta Math.},    volume={148},    date={1982},    pages={47--157},}




\bib{[L91]}{book}{
    author={Lawlor, G.R.},
    title={A sufficient criterion for a cone to be area-minimizing,},
   place={Mem. of the Amer. Math. Soc.},
   volume={91},
   date={1991},
}


\bib{BL}{article}{
    author={{Lawson, Jr.}, H.B.},
    title={The equivariant Plateau problem and interior regularity},
    journal={Trans. Amer. Math. Soc.},
    volume={173},
    date={1972},
    pages={231--249},
}

\bib{[LO77]}{article}{

    author={{Lawson, Jr.}, H.B.},
    author={Osserman, R.},
    title={Non-existence, non-uniqueness and irregularity of solutions to the minimal surface system},
    journal={Acta Math.},
    volume={139},
    date={1977},
    pages={1--17},
}



\bib{[M80]}{article}{

    author={Mashimo, K.},
    title={Degree of the standard isometric minimal immersions of complex projective spaces into spheres},
    journal={Tsukuba J. Math.},
    volume={4},
    date={1980},
    pages={133--145},
}

\bib{[M81]}{article}{

    author={Mashimo, K.},
    title={Degree of the standard isometric minimal immersions of the symmetric spaces of rank one into spheres},
    journal={Tsukuba J. Math.},
    volume={5},
    date={1981},
    pages={291--297},
}

\bib{O}{book}{
  author={Ohnita, Y.},
    title={The first standard minimal immersions of compact irreducible symmetric spaces},
   place={pp. 37--49,
   Differential geometry of submanifolds, 
   ed. by K. Kenmotsu, Lecture Notes in Mathematics, Vol. 1090, Springer-Verlag, 
},
                 date={1984},
}


\bib{PS1}{article}{
  author={Simoes, P.},
    title={On a class of minimal cones in $\mathbb R^n$},
    journal={Bull. Amer. Math. Soc.},
                 volume={3},
                 date={1974},
         pages={488--489}
}



\bib{PS2}{book}{
 author={Simoes, P.},
    title={A class of minimal cones in $\mathbb R^n$, $n\geq 8$, that minimize area},
    place={Ph.D. thesis, University of California, Berkeley, Calif.},
    date={1973},
   }

 \bib{LS}{book}{
    author={Simon, L.},
    title={Lectures on Geometric Measure Theory},
   place={Proc. Centre Math. Anal. Austral. Nat. Univ., Vol. 3,
   Centre for Mathematical Analysis, Canberra},
   date={1983},
}


 \bib{LS89}{article}{
    author={Simon, L.},
    title={Entire solutions of the minimal surface equation},
   journal={J. Diff. Geom.},
                 volume={30},
                 date={1989},
         pages={643--688}
}

\bib{Simons}{article}{
    author={Simons, J.},
    title={Minimal varieties in riemannian manifolds},
    journal={Ann. Math.},
    volume={88},
    date={1968},
    pages={62--105},
}


\bib{[T66]}{article}{
    author={Takahashi, T.},
    title={Minimal immersions of Riemannian manifolds},
    journal={J. Math. Soc. Japan},
    volume={18},
    date={1966},
    pages={380--385},
}

\bib{[T01]}{article}{
    author={Tang, Z.Z.},
    title={Nonexistence of a submersion from the 23-sphere to the Cayley projective plane},
    journal={Bull. London Math. Soc.},
    volume={33},
    date={2001},
    pages={347--350},
}



\bib{[U85]}{article}{
    author={Urakawa, H.},
    title={Minimal immersions of projective spaces into spheres},
    journal={Tsukuba J. Math.},
    volume={9},
    date={1985},
    pages={321--347},
}

\bib{W}{book}{
  author={Wallach, N.},
    title={Minimal immersions of symmetric spaces into spheres},
   place={pp. 1--40, Symmetric Spaces, Pure and Applied Mathematics, Vol. 8. Dekker, New York},
                 date={1972},
}


\bib{[XYZ16]}{article}{
author={Xu, X.W.}
author={Yang, L.}
   author={Zhang, Y.S.},
   title={On Lawson-Osserman constructions}
   journal={arXiv:1610.08162}
}
\end{biblist}
\end{bibdiv}
\end{document}